\documentclass[11pt]{article}       % onecolumn (second format)
\usepackage{amsmath}
\usepackage{amscd}
\usepackage{amssymb}
\usepackage{rotating}
\usepackage{amsfonts}
\usepackage{amsthm}
\usepackage{stmaryrd}
\usepackage[top=1in, bottom=1in, left=1.25in, right=1.25in]{geometry}

\input xy
\xyoption{all}

\allowdisplaybreaks[1]
\newcommand{\pd}{\partial}

\newcommand{\lbr}{{\llbracket}}
\newcommand{\rbr}{{\rrbracket}}

\DeclareMathOperator{\ch}{ch}

\newtheorem{theorem}{Theorem}[section]
\newtheorem{theorem/definition}{Theorem/Definition}[section]
\newtheorem{proposition}[theorem]{Proposition}
\newtheorem{lemma}[theorem]{Lemma}
\newtheorem{corollary}[theorem]{Corollary}
\newtheorem{Conjecture}{Conjecture}

\theoremstyle{remark}
\newtheorem{remark}{Remark}[section]

\theoremstyle{definition}
\newtheorem{definition}[theorem]{Definition}

\theoremstyle{definition}
 \newtheorem{example}{Example}[section]
\newcommand{\be}{\begin{equation}}
\newcommand{\ee}{\end{equation}}
\newcommand{\bea}{\begin{eqnarray}}
\newcommand{\eea}{\end{eqnarray}}
\newcommand{\ben}{\begin{eqnarray*}}
\newcommand{\een}{\end{eqnarray*}}
\newcommand{\bet}{\begin{equation}
\begin{split}}
\newcommand{\eet}{\end{split}
\end{equation}}

\usepackage{titlesec}

\begin{document}
\title
{\Large The Quantum McKay Correspondence for Singularities of Type D
\thanks{Date: July 19, 2012}}
\author{\normalsize Xiaowen Hu}
\date{}
\maketitle

\begin{abstract}
We prove the quantum McKay correspondence formulae conjectured by J.Bryan and A.Gholampour for the type D (binary) polyhedral groups in $SU(2)$ and $SO(3)$. We use
the method of induction by the WDVV equation and from the normal subgroups by J.Bryan and A.Gholampour , and the polynomiality technique developed in this article. We are also based on the validity of the corresponding conjecture for type A groups, which is proved by T.  Coates,  A. Corti,  H. Iritani,  and H.-H.  Tseng.\\

\end{abstract}

\section{Introduction}
Let $\mathcal{X}$ be a an effective orbifold with the coarse moduli space $X$ and $Y\rightarrow X$ be a crepant resolution. The general principle of McKay correspondence expects that the geometry of $\mathcal{X}$ coincides with that of $Y$. For the case of $[\mathbb{C}^2/G]$, where $G$ is a finite subgroup of $SU(2)$, the classical McKay correspondence states that the representation theory of $G$ is related to the intersection matrix of the exceptional divisors in $Y$. In the language of orbifold cohomology \cite{ChenRuan} \cite{Ruan}, this means a correspondence between the orbifold cohomology of $[\mathbb{C}^2/G]$ and the ordinary cohomology of $\widehat{\mathbb{C}^2/G}$. This correspondence is extended to a \emph{quantum} version, see \cite{BGr}, \cite{CoatesRuan}. We refer the reader also to see the introduction of \cite{Coates}, \cite{Zhou0811} for an ample discussion. When $\mathcal{X}$ satisfies the Hard Lefschetz condition, which is the case for the quotients of  binary polyhedral groups and polyhedral groups, it is conjectured \cite{BGr} that the (equivariant) primary orbifold Gromov-Witten invariants of $\mathcal{X}$ and the (equivariant) primary Gromov-Witten invariants of $Y$ is equal after change of variables and analytic continuation. \\

J. Bryan and A. Gholampour have computed the genus zero equivariant Gromov-Witten of $\widehat{\mathbb{C}^2/G}$ and $\widehat{\mathbb{C}^3/G}$ explicitly in \cite{BGh0707}, \cite{BGh0803} for $G\subset SU(2)$ and $G\subset SO(3)$, thus obtained a prediction for the  genus zero equivariant orbifold Gromov-Witten of $[\mathbb{C}^2/G]$. For their conjectural formulae, see Conjecture \ref{27} and Conjecture \ref{28}. Their conjecture is proved for the binary $E_6$ and $D_4$ groups by themselves \cite{BGh0708}, and by T. Coates, A. Corti, H. Iritani, and Hsian-Hua Tseng for the binary polyhedral and polyhedral groups of type A. The higher genus correspondence is partly proved by the results of \cite{Maulik} and \cite{Zhou0811}. \\

To prove J. Bryan and A. Gholampour's conjecture, one needs to compute some Hurwitz-Hodge integrals. For the case of non-abelian groups, the \emph{quantum Lefschetz principle} \cite{CG} \cite{Ts} does not applies, at least not directly. The method in \cite{BGh0708} is using the WDVV equation, and the reduction-to-normal subgroups method, to reduce the concerning integrals finally to the integrals computed in \cite{BGP}. We use this approach
to reduce the concerning integrals of binary type D groups to the integrals of type A, together with a small \emph{exceptional} class of integrals. Thus we can apply the result of \cite{CCIT}, \cite{CCIT2} to show that the conjecture holds for a infinite series of binary $D_n$. For computing the exceptional class of integrals, one can take the approach of \cite{Zhou0710}, which is based on the \emph{quantum Riemann-Roch theorem} of \cite{CG}, \cite{Ts}. But in general, this approach does not give a closed formula when one cannot control the number of Chern characters involved. This is partly because we don't have a prediction for the closed formula of \emph{descendent} Gromov-Witten invariants of
$[\mathbb{C}^2/G]$. In general, as in the case of manifolds, such closed formulae do not exist. One can only search for an \emph{integrable hierarchy} satisfied by the total descendent (or ancestor) Gromov-Witten potential. Such a prediction does not exist in our case up to now. For the case of $[\mathbb{C}^3/\mathbb{Z}_2]$, the reader may refer  \cite{Brini}, \cite{BCR}.\\

However, since the prediction for the exceptional integrals is a polynomial of $n-2$, one can try to prove this polynomiality,  and since we have proved the conjecture holds for infinitely many $n$, it holds \emph{for all} $n$. \\

The proof is organized as follows. In section 2 we recall the orbifold quantum Riemann-Roch theorem of \cite{CG}, \cite{Ts}, and give an explicit graph presentation of the formula at least for the genus 0 part, as an enhancement of the approach of \cite{Zhou0710}. In section 3 we recall the content of classical and quantum McKay correspondences, the conjectural formulae of Bryan and Gholampour \cite{BGh0707} and the result for the type A case \cite{CCIT2}. In section 4 we adopt the method of \cite{BGh0708} and show that the values of some special correlatos can be deduced from the corresponding formula of the normal groups. We call this procedure \emph{induction from normal subgroups}. In section 5 we adopt another method of \cite{BGh0708} and observe that all the correlators can be determined by some special correlators and the WDVV equation. We call this procedure \emph{induction by the WDVV equation}. Roughly speaking, we inductively solve the linear equations formed by the highest linear terms given by the equaliy of the coefficients of the WDVV equation. This together with the results of section 4 is enough to show that the quantum McKay correspondence holds for a infinite series of $\hat{D}_{n}$, namely for $n=2^{m}+2$, $m\geq 1$. At the end of this section, we sketch a proof of the conjecture for $D_n\subset SO(3)$. In section 6 we associate a number $S_{n}(\Gamma)$ to every \emph{decorated Feymann diagram}, and prove the (Laurent) polynomial dependence on $n$. In section 7 we rearrange the form of the quantization operator in the orbifold quantum Riemann-Roch in the new coordinates, then we apply the graph presentation of section 1  and further reduce it to \emph{a summation of fractional Bernoulli numbers over graphs} which is treated in section 6. Thus the result of section 6 gives the polynomiality we need. \\

One may hope to take a similar approach to obtain the quantum McKay correspondence of type D for the \emph{full} gravitational genus zero correlators.  For this, one may need to firstly get a conjectural form of the $J$-functions, for which the quantum Lefschetz doesn't directly apply. Then one can try to do the normal group induction and the WDVV (replaced by the genus zero topological recursion relations). The polynomiality theorem still holds, as remarked at the end of section 7. Since our approach needs the genus to be zero in many aspects, the cases of higher genera remain open.\\

\textit{Acknowledgements.} The author thanks Prof.~Jian Zhou for his great patience and guidance during all the time. He also thanks Huazhong Ke, Di Yang for helpful discussions. The author is particularly indebted to Hanxiong Zhang, who gives him kind encouragement and useful suggestions through the writing.

\section{Equivariant Gromov-Witten invariants of $[\mathbb{C}^n /G]$ and orbifold quantum Riemann-Roch }
\subsection{Definition of the Equivariant Gromov-Witten invariants }
Let a finite group $G$ act on $\mathbb{C}^n$ via a representation $\rho$ of $G$, and let $\mathbb{C}^{*}$ act on $\mathbb{C}^n$ diagonally. When we want to emphasize that $\mathbb{C}^n$ is the representation space of $\rho$, we denote it by $V_{\rho}$.
The action of $\mathbb{C}^{*}$ commutes with the action of $G$, therefore descends to be an action on $[\mathbb{C}^n /G]$, and has a unique isolated fixed substack $\mathcal{B}G$. We define the equivariant Gromov-Witten invariants of $[\mathbb{C}^n /G]$ by virtual localization \cite{GR}. When the monodromy data $\alpha$ is not trivial at every marked points, the moduli space $\overline{\mathcal{M}}_{g,\alpha}([\mathbb{C}^n /G])$ is proper, and the definition coincides with the usual definition \cite{AGV}. For details, see, e.g., \cite{Zhou0811}.\\
We denote the equivariant parameter by $\lambda$, the equivariant Euler class by $c_{T}(\cdot)$.  We use the notations in \cite{Zhou0710}, \cite{Zhou0811} and using virtual localization \cite{GR}, the equivariant correlators are given by
\bea\label{17}
\langle \prod_{j=1}^{m}\tau_{k_j}(e_{\lbr \gamma_{j}\rbr})\rangle^{[\mathbb{C}^n /G]}=\int_{\overline{\mathcal{M}}_{g,m}(\mathcal{B}G;\coprod_{i=1}^{m}\lbr \gamma_{i}\rbr)}
\frac{c_{T}(\mathbb{F}_{\rho,g,m}^{1})}{c_{T}(\mathbb{F}_{\rho,g,m}^{0})}\cdot \prod_{j=1}^{m}\bar{\psi}_{j}^{k_j}.
\eea

Another way to define equivariant Gromov-Witten invariants of $[\mathbb{C}^n /G]$ is  through Givental's formalism of \textit{twisted} Gromov-Witten invariants, see \cite{CG}, \cite{Ts}. Fix a multiplicative characteristic class $\mathbf{c}(\cdot)=\exp (\sum_{k=0}^{\infty}s_{k}\ch_{k}(\cdot))$, where $s_{0},s_{1},\cdots$ are parameters. Following \cite{Ts}, define the $(\mathbf{c},\rho)$-\textit{twisted Gromov-Witten total descendent potential} of $\mathcal{B}G$ by
\bea\label{18}
Z^{\mathbf{c},\rho}(\mathbf{t})= \exp\Bigg(\sum_{g=0}^{\infty}\hbar^{2g-2}\sum_{n\geq 0}\frac{1}{n!}
\int_{\overline{\mathcal{M}}_{g,n}(\mathcal{B}G)}\mathbf{c}(\mathbb{F}_{\rho,g,n})\prod_{i=1}^{n}\sum_{k=1}^{\infty}ev_{i}^{*}(\sum_{\lbr \gamma \rbr}t_{k}^{\lbr \gamma\rbr}e_{\lbr\gamma\rbr})\bar{\psi}_{i}^{k}\Bigg).
\eea
Then  specializing the parameters to $s_{0}=\ln \lambda$, and $s_{k}=-(k-1)!/\lambda^{k}$ for $k\geq 1$, we obtain the equivariant Gromov-Witten invariants. \\

The two definitions give the same invariants, since the relation between \emph{equivariant} Chern classes and Chern characters (i.e., the \textit{Newton's identities}) holds for \textit{virtual} vector bundles.

\subsection{The orbifold quantum Riemann-Roch theorem}
Givental's ingenious \textit{quantization formalism} makes a way to relate the twisted Gromov-Witten invariants to the ordinary ones. We recall the formula in the following, and refer the readers to \cite{CG}, \cite{Ts}, \cite{Zhou0710} for the notations. Note that the coefficient of $s_{0}=\ln \lambda$ encodes the virtual \textit{rank} of the virtual \textit{Hurwitz-Hodge} bundle, and in the practical computations this rank on the component which we concern is known, thus we can set $s_{0}=0$ and multiply appropriate powers of $\lambda$ to the resulting integrals. In this way, the quantum Riemann-Roch formula becomes more simple.
\begin{theorem}[\cite{CG}, \cite{Ts}]\label{18}
For $p\geq 1$, we have
\ben
\frac{\pd Z^{\mathbf{c},\rho}(\mathbf{t})}{\pd s_p}= \Bigg(s_{p}\frac{A_{p+1}(V_{\rho})z^p}{(p+1)!}\Bigg)^{\wedge}Z^{\mathbf{c},\rho}.\\
\een
\end{theorem}
This formula provides an effective way to compute Hurwitz-Hodge integrals, as is shown by J. Zhou in \cite{Zhou0710}. We recall this method. Denote by $F^{\mathbf{c},\rho}(\mathbf{t})$  the free energy
\ben
F^{\mathbf{c},\rho}(\mathbf{t})&=&\sum_{g=0}^{\infty}\hbar^{2g-2}F_{g}^{\mathbf{c},\rho}(\mathbf{t})\\
&=&\sum_{g=0}^{\infty}\hbar^{2g-2}\sum_{n\geq 0}\frac{1}{n!}
\int_{\overline{\mathcal{M}}_{g,n}(\mathcal{B}G)}\mathbf{c}(\mathbb{F}_{\rho,g,n})\prod_{i=1}^{n}\sum_{k=1}^{\infty}ev_{i}^{*}(\sum_{\lbr \gamma \rbr}t_{k}^{\lbr \gamma\rbr}e_{\lbr\gamma\rbr})\bar{\psi}_{i}^{k}.
\een
In the same way, we denote by $F^{G}(\mathbf{t})=\sum_{g=0}^{\infty}\hbar^{2g-2}F_{g}^{G}(\mathbf{t})$  the untwisted free energy.
The operator $\mathcal{O}_p (\rho):=\Big(\frac{A_{p+1}(V_{\rho})z^p}{(p+1)!}\Big)^{\wedge}$ is of the form
\bea\label{37}
\mathcal{O}_p (\rho)=D_{p} +\frac{\hbar^2}{2}\sum_{l=0}^{p-1}C_{p}^{\alpha\beta}\pd_{\beta, l}\pd_{\alpha,p-1-l},
\eea
where $D_p$ is a first order differential operator. Therefore
\bea
&&(Z^{\mathbf{c},\rho})^{-1}\Bigg(\frac{A_{p+1}(V_{\rho})z^p}{(p+1)!}\Bigg)^{\wedge}Z^{\mathbf{c},\rho}\nonumber\\
&=& (Z^{\mathbf{c},\rho})^{-1}(D_{p} +\frac{\hbar^2}{2}\sum_{l=0}^{p-1}C_{p}^{\alpha\beta}\pd_{\beta, l}\pd_{\alpha,p-1-l})\exp F^{\mathbf{c},\rho}\nonumber\\
&=& D_p F^{\mathbf{c},\rho} +\frac{\hbar^2}{2}\sum_{l=0}^{p-1}C_{p}^{\alpha\beta}\pd_{\beta, l}\pd_{\alpha,p-1-l}F^{\mathbf{c},\rho}+\frac{\hbar^2}{2}\sum_{l=0}^{p-1}C_{p}^{\alpha\beta}\pd_{\beta, l}F^{\mathbf{c},\rho}\pd_{\alpha,p-1-l}F^{\mathbf{c},\rho}.
\eea
Taking the coefficients of $\hbar^{-2}$, we get
\bea\label{24}
\frac{\pd}{\pd s_p}F_{0}^{\mathbf{c},\rho}&=&\sum_{n\geq 3}\frac{1}{n!}
\int_{\overline{\mathcal{M}}_{0,n}(\mathcal{B}G)}\ch_{p}(\mathbb{F}_{\rho,0,n})\mathbf{c}(\mathbb{F}_{\rho,0,n})\prod_{i=1}^{n}\sum_{k=1}^{\infty}ev_{i}^{*}(\sum_{\lbr \gamma \rbr}t_{k}^{\lbr \gamma\rbr}e_{\lbr\gamma\rbr})\bar{\psi}_{i}^{k}\nonumber\\
&=& D_p F_{0}^{\mathbf{c},\rho} +\frac{1}{2}\sum_{l=0}^{p-1}C_{p}^{\alpha\beta}\pd_{\beta, l}F_{0}^{\mathbf{c},\rho}\pd_{\alpha,p-1-l}F_{0}^{\mathbf{c},\rho}.
\eea
Then taking all $s_k=0$ we obtain
\bea
&&\sum_{n\geq 3}\frac{1}{n!}\int_{\overline{\mathcal{M}}_{0,n}(\mathcal{B}G)}\ch_{p}\prod_{i=1}^{n}\sum_{k=1}^{\infty}ev_{i}^{*}(\sum_{\lbr \gamma \rbr}t_{k}^{\lbr \gamma\rbr}e_{\lbr\gamma\rbr})\bar{\psi}_{i}^{k}\nonumber\\
&=&D_p F_{0}^{G} +\frac{1}{2}\sum_{l=0}^{p-1}C_{p}^{\alpha\beta}\pd_{\beta, l}F_{0}^{G}\pd_{\alpha,p-1-l}F_{0}^{G}.
\eea
We use the diagram
$$ 
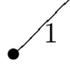
\xy
(0,0); (8,8), **@{-};(5,3)*+{1};
(0,0)*+{\bullet};
\endxy
$$
to indicate the first  term $D_p F_{0}^{G}$, and use the diagram
\bea\label{38}
\xy
(0,0); (10,0), **@{-};(5,2)*+{1};
(0,0)*+{\bullet};(10,0)*+{\bullet};
\endxy
\eea
to indicate the second term $\frac{1}{2}\sum_{l=0}^{p-1}C_{p}^{\alpha\beta}\pd_{\beta, l}F_{0}^{G}\pd_{\alpha,p-1-l}F_{0}^{G}$.
We give an example.
\begin{example}\label{19}
Consider $[\mathbb{C}^2 /\hat{D}_n]$ for an even $n\geq 4$. See section 3 and 4 for the definition and some coefficients in the following computations. By definition of the quantization of a symplectic transform,
\ben
&&\big(\frac{A_{p+1}(V_{\rho_{1}})z^p}{(p+1)!}\big)^{\wedge}\\
&=&\frac{2B_{p+1}}{(p+1)!}\pd_{\lbr 1\rbr,1+p}
-\frac{2B_{p+1}}{(p+1)!}\sum_{l= 0}^{\infty}t_{l}^{\lbr 1\rbr}\pd_{\lbr 1\rbr,l+p}-\sum_{k=1}^{n-3}\frac{B_{p+1}(\frac{k}{2n-4})+B_{p+1}(\frac{2n-4-k}{2n-4})}{(p+1)!}\sum_{l= 0}^{\infty}t_{l}^{\lbr a^k\rbr}\pd_{\lbr a^k\rbr,l+p}\\
&&-\frac{2B_{p+1}(\frac{n-2}{2n-4})}{(p+1)!}\sum_{l=0}^{\infty}t_{l}^{\lbr a^{n-2}\rbr}\pd_{\lbr a^{n-2}\rbr,l+p}-\frac{B_{p+1}(\frac{1}{4})+B_{p+1}(\frac{3}{4})}{(p+1)!}\sum_{l=0}^{\infty}(t_{l}^{\lbr b\rbr}\pd_{\lbr b\rbr,l+p}+t_{l}^{\lbr ab\rbr}\pd_{\lbr ab\rbr,l+p})
\\
&&+\frac{\hbar ^2}{2}\sum_{l= 0}^{p-1}(-1)^l\Bigg((4n-8)\frac{2B_{p+1}}{(p+1)!}\pd_{\lbr 1\rbr,l}\pd_{\lbr 1\rbr,p-1-l}+(2n-4)\sum_{k=1}^{n-3}\frac{B_{p+1}(\frac{k}{2n-4})+B_{p+1}(\frac{2n-4-k}{2n-4})}{(p+1)!}\\
&&\cdot\pd_{\lbr a^k \rbr,l}\pd_{\lbr a^k\rbr,p-1-l}
\\&&
+(4n-8)\frac{2B_{p+1}(\frac{1}{2})}{(p+1)!}\pd_{\lbr a^{n-2}\rbr,l}\pd_{\lbr a^{n-2}\rbr,p-1-l}+
4\cdot\frac{B_{p+1}(\frac{1}{4})+B_{p+1}(\frac{3}{4})}{(p+1)!}(\pd_{\lbr b\rbr,l}\pd_{
\lbr b \rbr,p-1-l}\\
&&+\pd_{\lbr ab\rbr,l}\pd_{
\lbr ab \rbr,p-1-l})\Bigg).\\
\een
Take $s_0=s_{1}=\cdots=0$ in (\ref{24}), we obtain
\ben
&&\int_{\overline{M}_{0,2m}(\mathcal{B}\hat{D}_{n};\lbr b\rbr^{2m})}\ch_{2m-3}(\mathbb{F}_{\rho,0,2m})\\
&=&\frac{2B_{2m-2}}{(2m-2)!}\int_{\overline{M}_{0,2m+1}(\mathcal{B}\hat{D}_{n};\lbr b\rbr^{2m},\lbr 1\rbr)}\bar{\psi}_{2m+1}^{2m-2}-\frac{B_{2m-2}(\frac{1}{4})+B_{2m-2}(\frac{3}{4})}{(2m-2)!}\cdot 2m\int_{\overline{M}_{0,2m}(\mathcal{B}\hat{D}_{n};\lbr b\rbr^{2m})}\bar{\psi}_{1}^{2m-3}\\
&&+ \frac{1}{2}\cdot(4n-8)\frac{2B_{2m-2}}{(2m-2)!}\sum_{i=1}^{m-1}\binom{2m}{2i}\int_{\overline{M}_{0,2i+1}(\mathcal{B}\hat{D}_{n};\lbr b\rbr^{2i},\lbr 1\rbr)}\bar{\psi}_{2i+1}^{2i-2}\int_{\overline{M}_{0,2m-2i+1}(\mathcal{B}\hat{D}_{n};\lbr b\rbr^{2m-2i},\lbr 1\rbr)}\bar{\psi}_{2m-2i+1}^{2m-2i-2}\\
&&+ \frac{1}{2}\cdot(2n-4)\sum_{k=1}^{n-3}\frac{B_{2m-2}(\frac{k}{2n-4})+B_{2m-2}(\frac{2n-4-k}{2n-4})}{(2m-2)!}
\sum_{i=1}^{m-1}\binom{2m}{2i}\int_{\overline{M}_{0,2i+1}(\mathcal{B}\hat{D}_{n};\lbr b\rbr^{2i},\lbr a^{k}\rbr)}\bar{\psi}_{2i+1}^{2i-2}\\
&&\cdot\int_{\overline{M}_{0,2m-2i+1}(\mathcal{B}\hat{D}_{n};\lbr b\rbr^{2m-2i},\lbr a^{k}\rbr)}\bar{\psi}_{2m-2i+1}^{2m-2i-2}+\frac{1}{2}\cdot(4n-8)\frac{2B_{2m-2}(\frac{1}{2})}{(2m-2)!}\sum_{i=1}^{m-1}\binom{2m}{2i}\\
&&\cdot\int_{\overline{M}_{0,2i+1}
(\mathcal{B}\hat{D}_{n};\lbr b\rbr^{2i},\lbr a^{n-2}\rbr)}\bar{\psi}_{2i+1}^{2i-2}\int_{\overline{M}_{0,2m-2i+1}(\mathcal{B}\hat{D}_{n};\lbr b\rbr^{2m-2i},\lbr a^{n-2}\rbr)}\bar{\psi}_{2m-2i+1}^{2m-2i-2}\\
&&+2\cdot\frac{B_{2m-2}(\frac{1}{4})+B_{2m-2}(\frac{3}{4})}
{(2m-2)!}\sum_{i=2}^{2m-2}(-1)^{i}\binom{2m}{i}\int_{\overline{M}_{0,i+1}
(\mathcal{B}\hat{D}_{n};\lbr b\rbr^{i},\lbr b\rbr)}\bar{\psi}_{i+1}^{i-2}\\
&&\cdot\int_{\overline{M}_{0,2m-i+1}(\mathcal{B}\hat{D}_{n};\lbr b\rbr^{2m-i},\lbr b\rbr)}\bar{\psi}_{2m-i+1}^{2m-i-2}\\
&=&\frac{2B_{2m-2}}{(2m-2)!}\frac{(n-2)^{2m-1}}{4n-8}-\frac{B_{2m-2}(\frac{1}{4})+B_{2m-2}(\frac{3}{4})}{(2m-2)!}\cdot 2m\frac{(n-2)^{2m-1}}{4n-8}\\
&& +(2n-4)\cdot\frac{2B_{2m-2}}{(2m-2)!}\sum_{i=1}^{m-1}\binom{2m}{2i}\frac{(n-2)^{2i-1}}{4n-8}\frac{(n-2)^{2m-2i-1}}{4n-8}\\
&&+(n-2)\sum_{j=1}^{\frac{n}{2}-2}\frac{B_{2m-2}(\frac{2j}{2n-4})+B_{2m-2}(\frac{2n-4-2j}{2n-4})}{(2m-2)!}
\sum_{i=1}^{m-1}\binom{2m}{2i}\frac{2(n-2)^{2i-1}}{4n-8}\frac{2(n-2)^{2m-2i-1}}{4n-8}\\
&&+(2n-4)\cdot\frac{2B_{2m-2}(\frac{1}{2})}{(2m-2)!}\sum_{i=1}^{m-1}\binom{2m}{2i}\frac{(n-2)^{2i-1}}{4n-8}\frac{(n-2)^{2m-2i-1}}{4n-8}\\
&&+2\cdot\frac{B_{2m-2}(\frac{1}{4})+B_{2m-2}(\frac{3}{4})}{(2m-2)!}\sum_{j=1}^{m-2}(-1)\binom{2m}{2j+1}
\frac{(n-2)^{2j+1}}{4n-8}\frac{(n-2)^{2m-2j-1}}{4n-8}\\
&=&\frac{B_{2m-2}}{(2m-2)!}\frac{(n-2)^{2m-2}}{2}-
\frac{(n-2)^{2m-2}}{8}\cdot\frac{B_{2m-2}(\frac{1}{4})+B_{2m-2}(\frac{3}{4})}{(2m-2)!}\sum_{j=0}^{m-1}\binom{2m}{2j+1}\\
&&+\frac{(n-2)^{2m-3}}{4}\sum_{j=0}^{n-3}\frac{B_{2m-2}(\frac{2j}{2n-4})}{(2m-2)!}\sum_{i=1}^{m-1}\binom{2m}{2i}\\
&=&\frac{B_{2m-2}}{(2m-2)!}\Big((n-2)^{2m-2}+2^{2m-3}\Big)\Big(1-\frac{1}{2^{2m-2}}\Big).\\
\een
In the last equality we used
\ben
\sum_{k=0}^{n-1}B_p (\frac{k}{n})=\frac{1}{n^{p-1}}B_p,\\
\sum_{k=0}^{n-1}B_p (\frac{2k+1}{2n})=\Big(\frac{1}{(2n)^{p-1}}-\frac{1}{n^{p-1}}\Big)B_p.\\
\een
When $m=4$, this value coincides with the conjecture \ref{27}, see section 4.2.
\hfill\qedsymbol
\end{example}

For the integrals involving more than one Chern characters, we can further differentiate (\ref{24}) by the $s_k$'s, then set all $s_k=0$. For example,
\ben
&&\frac{\pd}{\pd s_{p_2}}\frac{\pd}{\pd s_{p_1}}F_{0}^{\mathbf{c},\rho}\\&=&D_{p_1}\frac{\pd}{\pd s_{p_2}} F_{0}^{\mathbf{c},\rho} +\frac{1}{2}\sum_{l=0}^{p_{1}-1}C_{p_1}^{\alpha\beta}\pd_{\beta, l}\Big(\frac{\pd}{\pd s_{p_2}}F_{0}^{\mathbf{c},\rho}\Big)\pd_{\alpha,p_{1}-1-l}F_{0}^{\mathbf{c},\rho}\\
&&+\frac{1}{2}\sum_{l=0}^{p_{1}-1}C_{p_1}^{\alpha\beta}\pd_{\beta, l}F_{0}^{\mathbf{c},\rho}\pd_{\alpha,p_{1}-1-l}\Big(\frac{\pd}{\pd s_{p_2}}F_{0}^{\mathbf{c},\rho}\Big)\\
&=&D_{p_1}\Big(D_{p_2} F_{0}^{\mathbf{c},\rho} +\frac{1}{2}\sum_{l=0}^{p_2-1}C_{p_2}^{\alpha_2\beta_2}\pd_{\beta_2, l}F_{0}^{\mathbf{c},\rho}\pd_{\alpha_2,p_2-1-l}F_{0}^{\mathbf{c},\rho}\Big) +\frac{1}{2}\sum_{l=0}^{p_{1}-1}C_{p_1}^{\alpha_1\beta_1}\pd_{\beta_1, l}\Big(D_{p_2} F_{0}^{\mathbf{c},\rho} \\
&&+\frac{1}{2}\sum_{l=0}^{p_2-1}C_{p_2}^{\alpha_2\beta_2}\pd_{\beta_2, l}F_{0}^{\mathbf{c},\rho}\pd_{\alpha_2,p-1-l}F_{0}^{\mathbf{c},\rho}\Big)\pd_{\alpha_1,p_{1}-1-l}F_{0}^{\mathbf{c},\rho}\\
&&+\frac{1}{2}\sum_{l=0}^{p_{1}-1}C_{p_1}^{\alpha_1\beta_1}\pd_{\beta_1, l}F_{0}^{\mathbf{c},\rho}\pd_{\alpha_1,p_{1}-1-l}\Big(D_{p_2} F_{0}^{\mathbf{c},\rho} +\frac{1}{2}\sum_{l=0}^{p_2-1}C_{p_2}^{\alpha_2\beta_2}\pd_{\beta_2, l}F_{0}^{\mathbf{c},\rho}\pd_{\alpha_2,p_2-1-l}F_{0}^{\mathbf{c},\rho}\Big),
\een
taking all $s_k=0$ we obtain
\ben
&&\sum_{n\geq 3}\frac{1}{n!}\int_{\overline{\mathcal{M}}_{0,n}(\mathcal{B}G)}\ch_{p_1}\ch_{p_2}\prod_{i=1}^{n}\sum_{k=1}^{\infty}ev_{i}^{*}(\sum_{\lbr \gamma \rbr}t_{k}^{\lbr \gamma\rbr}e_{\lbr\gamma\rbr})\bar{\psi}_{i}^{k}\\
&=&(D_{p_1}\circ D_{p_2}) F_{0}^{G}+ \sum_{l=0}^{p_2-1}C_{p_2}^{\alpha_2\beta_2}(D_{p_1}\circ \pd_{\beta_2, l})F_{0}^{G}\cdot\pd_{\alpha_2,p_2-1-l}F_{0}^{G}\\
&&+\sum_{l=0}^{p_{1}-1}C_{p_1}^{\alpha_1\beta_1}(\pd_{\beta_1, l}\circ D_{p_2} )F_{0}^{G} \cdot\pd_{\alpha_1,p_{1}-1-l}F_{0}^{G}\\
&&+\sum_{l_1=0}^{p_{1}-1}\sum_{l_2=0}^{p_2-1}C_{p_1}^{\alpha_1\beta_1}C_{p_2}^{\alpha_2\beta_2}(\pd_{\beta_1, l_1}\pd_{\beta_2, l_2}F_{0}^{G})\pd_{\alpha_2,p-1-l_2}F_{0}^{G}\pd_{\alpha_1,p_{1}-1-l_1}F_{0}^{G}.
\een
We use the diagram
$$ 
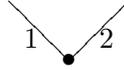
\xy
(0,0); (8,8), **@{-};(5,3)*+{2};
(0,0); (-8,8), **@{-};(-5,3)*+{1};
(0,0)*+{\bullet};
\endxy
$$
to indicate the first  term $(D_{p_1}\circ D_{p_2}) F_{0}^{G}$, and use the diagram
$$ 
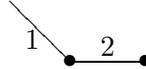
\xy
(0,0); (10,0), **@{-};(5,2)*+{2};
(0,0); (-8,8), **@{-};(-5,3)*+{1};
(0,0)*+{\bullet};(10,0)*+{\bullet};
\endxy
$$
to indicate the second term $\sum_{l=0}^{p_2-1}C_{p_2}^{\alpha_2\beta_2}(D_{p_1}\circ \pd_{\beta_2, l})F_{0}^{G}\cdot\pd_{\alpha_2,p_2-1-l}F_{0}^{G}$,  and use the diagram
$$ 
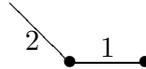
\xy
(0,0); (10,0), **@{-};(5,2)*+{1};
(0,0); (-8,8), **@{-};(-5,3)*+{2};
(0,0)*+{\bullet};(10,0)*+{\bullet};
\endxy
$$
to indicate the third term $\sum_{l=0}^{p_{1}-1}C_{p_1}^{\alpha_1\beta_1}(\pd_{\beta_1, l}\circ D_{p_2} )F_{0}^{G} \cdot\pd_{\alpha_1,p_{1}-1-l}F_{0}^{G}$, and use the diagram
$$ \xy
(0,0); (10,0), **@{-};(5,3)*+{1};
(10,0); (20,0), **@{-};(15,3)*+{2};
(0,0)*+{\bullet};(10,0)*+{\bullet};(20,0)*+{\bullet};
\endxy
$$
to indicate the fourth term $\sum_{l_1=0}^{p_{1}-1}\sum_{l_2=0}^{p_2-1}C_{p_1}^{\alpha_1\beta_1}C_{p_2}^{\alpha_2\beta_2}(\pd_{\beta_1, l_1}\pd_{\beta_2, l_2}F_{0}^{G})\pd_{\alpha_2,p-1-l_2}F_{0}^{G}\pd_{\alpha_1,p_{1}-1-l_1}F_{0}^{G}$.

\begin{definition}
A \emph{labeled tree with half edges} (abbreviated by \emph{LTHE} ) is a connected tree $\Gamma$ with edges and half edges, together with a \emph{bijection} from the set of edges and half edges to the set $\{1,\cdots, n\}$. This bijection is called the \emph{labeling}, and $n$ called the \emph{degree} of the labeled tree with half edges.
\end{definition}
Let $\{\mathcal{O}_{p}\}_{1\leq p\leq r}$ be an ordered set of differential operators of the form (\ref{37}), where $C_{p}^{\alpha\beta}$ are constant. For a function of $\mathbf{t}$, $F(\mathbf{t})$, consider the expansion of
\ben
\mathcal{O}_{1}\cdots\mathcal{O}_{r}F(\mathbf{t}),
\een
one sees that each labeled tree with half edges  $\Gamma$ with degree $r$ gives a contribution, which we denote by $\text{Cont}(\Gamma;\mathcal{O}_{1}\cdots\mathcal{O}_{r}F(\mathbf{t}))$. Note that the only labeled tree with half edges that has a nontrivial
automorphism is the diagram (\ref{38}), whose automorphism group is $\mathbb{Z}_2$, and the factor of $\frac{1}{2}$ is included in the contribution.\\

\begin{remark}
The differential operators $\mathcal{O}_{p}(\rho)$'s coming from quantization commute pairwisely, while the first and second derivatives parts of $\mathcal{O}_{p}(\rho)$'s do not commute. But note that in our setting, once the order of $\mathcal{O}_{1},\cdots, \mathcal{O}_{r}$ is given, the order of the differentiation for every LTHE is determined with no ambiguity.
\end{remark}

We conclude that
\begin{theorem}\label{47}
\ben
&&\sum_{n\geq 3}\frac{1}{n!}\int_{\overline{\mathcal{M}}_{0,n}(\mathcal{B}G)}\ch_{p_1}\cdots\ch_{p_r}\prod_{i=1}^{n}\sum_{k=1}^{\infty}ev_{i}^{*}(\sum_{\lbr \gamma \rbr}t_{k}^{\lbr \gamma\rbr}e_{\lbr\gamma\rbr})\bar{\psi}_{i}^{k}\\
&=&\sum_{\emph{deg}(\Gamma)=r}\emph{Cont}(\Gamma; \mathcal{O}_{p_1}(\rho)\cdots\mathcal{O}_{p_r}(\rho)F_{0}^{G}),
\een
where $\Gamma$ runs over the set of labeled tree with half edges of degree $r$.
\end{theorem}
\hfill\qedsymbol

\section{McKay correspondence and Bryan-Gholampour conjecture}
\subsection{The classical McKay correspondence}
For a finite subgroup $G$ of $SU(2)$, J. McKay \cite{McKay} observed the connection between the representation theory of $G$ and the Dynkin diagram arising from the configuration of the exceptional divisors of  the crepant resolution of $\mathbb{C}^2/G$. Here we recall explicitly the McKay correspondence for the binary dihedral group $\hat{D}_{n}$. We give the details for the reader's convenience and to fix the notations. The readers can also refer to \cite{GSV}, \cite{Reid}.\\
The binary dihedral group $\hat{D}_n$ is generated by $a$, $b$, with relations $a^{n-2}=b^{2}$, $b^{4}=1$, $ba=a^{-1}b$. There are four 1-dimensional representations of $\hat{D}_n$. They are\\
\begin{center}\begin{tabular}{l|c|c}\hline
 & $a$ & $b$\\ \hline
$\psi_{1}$ &  $1$ & $1$ \\ \hline
$\psi_{2}$ &  1 & $-1$ \\ \hline
$\psi_{3}$ &  $-1$ & $-1$ \\ \hline
$\psi_{4}$ &  $-1$ & 1 \\ \hline
\end{tabular}\end{center}
for $2|n$, and \\
\begin{center}\begin{tabular}{l|c|c}\hline
 & $a$ & $b$\\ \hline
$\psi_{1}$ &  $1$ & $1$ \\ \hline
$\psi_{2}$ &  1 & $-1$ \\ \hline
$\psi_{3}$ &  $-1$ & $\sqrt{-1}$ \\ \hline
$\psi_{4}$ &  $-1$ & $-\sqrt{-1}$\\ \hline
\end{tabular}\end{center}
for $2\nmid n$. There are $n-3$ 2-dimensional representations $\rho_{1},\cdots,\rho_{n-3}$, given by
\bea
\rho_{k}(a)= \left (\begin{array}{cc}
\omega^{k} & \\
& \omega^{-k}
\end{array}
\right), &
\rho_{k}(b)= \left (\begin{array}{cc}
 & 1\\
(-1)^{k}&
\end{array}
\right),
\eea
where $\omega=\exp (\frac{2\pi \sqrt{-1}}{2n-4})$. These $n+1$ representations form the set of the complex irreducible representations of $\hat{D}_{n}$. The inclusion $\hat{D}_{n}\subset SU(2)$ is given by $\rho_{1}$. Consider the tensor product of $\rho_{1}$ with the other irreducible representations, we have
\bea
\rho_{1}\otimes \psi_1=\rho_{1}, & \rho_{1}\otimes \psi_2=\rho_{1},\\
\rho_{1}\otimes \psi_3=\rho_{n-3}, & \rho_{1}\otimes \psi_4=\rho_{n-3},\\
\rho_{1}\otimes \rho_{1}=\rho_{2}\oplus \psi_{1}\oplus\psi_{2}, &\rho_{1}\otimes \rho_{n-3}=\rho_{n-4}\oplus \psi_{3}\oplus\psi_{4},\\
\rho_{1}\otimes \rho_{k}= \rho_{k-1}\oplus \rho_{k+1},& \quad \text{for}\quad 2\leq k\leq n-4.\\
\eea
These decompositions correspond to the (extended) Dynkin diagram
$$ 
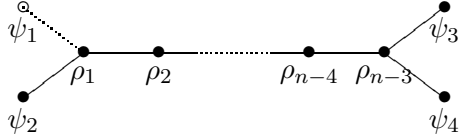
\xy
(0,0); (10,0), **@{-};(15,0), **@{-}; (25,0), **@{.}; (30,0),**@{-};(40, 0), **@{-};
(0,0)*+{\bullet};(10,0)*+{\bullet}; (30,0)*+{\bullet}; (40,0)*+{\bullet};
(0,-3)*+{\rho_1};(10,-3)*+{\rho_2};(30,-3)*+{\rho_{n-4}}; (40,-3)*+{\rho_{n-3}};
(0,0); (-8,6), **@{.};(0,0); (-8,-6), **@{-};(-8,-6)*+{\bullet};(-8,-9)*+{\psi_2};(-8,3)*+{\psi_1};(-8,6)*+{\circ};
(40,0); (48,6), **@{-};(40,0); (48,-6), **@{-};(48,-6)*+{\bullet};(48,-9)*+{\psi_4};(48,3)*+{\psi_3};(48,6)*+{\bullet};
\endxy
$$
Consider also the crepant resolution $Y^{\hat{D}_{n}}$ of $\mathbb{C}^2/\hat{D}_{n}$, we draw a node for each irreducible component of the exceptional divisor, and draw an edge connecting two nodes when the corresponding components intersects, then we obtain the same Dynkin diagram. More precisely, the intersection matrix is the minus Cartan matrix. This is the classical McKay correspondence.\\

\subsection{The quantum McKay correspondence}
As a special case of the general crepant resolution conjecture \cite{BGr}, \cite{CoatesRuan}, \cite{Ruan}, Bryan and Gholampour \cite{BGh0707} made the following
\begin{Conjecture}\label{27}
Let  $F_{0}^{\mathcal{X}}(x_1,\cdots,x_n)$  denote  the  $\mathbb{C}^{*} $-equivariant  genus 0 orbifold Gromov-Witten potential  of the orbifold $\mathcal{X}= [\mathbb{C}^{2}/G]$,  where we have set the unit parameter $x_{0}$  equal to zero. Let $R$ be the root system associated to G. Denote by $\alpha_{1},\cdots, \alpha_{n}$ the simple roots, and $R^{+}$ the set of positive roots. Then
\bea\label{25}
F_{0}^{\mathcal{X}}(x_1,\cdots,x_n)=2\lambda\sum_{\beta\in R^{+}}h(\pi+P_{\beta}),
\eea
where $h(u)$ is a series with
\bea
h^{\prime\prime\prime}(u)=\frac{1}{2}\tan(\frac{-u}{2})
\eea
and
\bea
P_{\beta}=\sum_{k=1}^{n}\frac{b_{k}}{|G|}\Bigg(2\pi n_{k}+\sum_{g\in G}\sqrt{2-\chi_{\rho_{1}}(g)}\overline{\chi}_{k}(g)x_{\lbr g\rbr}\Bigg)
\eea
where $b_k$   are  the coefficient of $\beta=\sum_{k}b_{k}\alpha_{k}$ and $n_k$   are  the coefficients of the largest  root. Note that $n_k$ is also the dimension of $\chi_{k}$.
\end{Conjecture}
For the polyhedral subgroups $G$ $\subset SO(3)$, Bryan and Gholampour \cite{BGh0803} made the following
\begin{Conjecture}\label{28}
Let  $F_{0}^{\mathcal{X}}(x_1,\cdots,x_n)$  denote  the  $\mathbb{C}^{*} $-equivariant  genus 0 orbifold Gromov-Witten potential  of the orbifold $\mathcal{X}= [\mathbb{C}^{3}/G]$,  where we have set the unit parameter $x_{0}$  equal to zero.  Then
\bea\label{26}
F_{0}^{\mathcal{X}}(x_1,\cdots,x_n)=\frac{1}{2}\sum_{\beta\in R^{+}}h(\pi+P_{\beta}),
\eea
where
\bea
P_{\beta}=\sum_{\rho}\frac{b_{\rho}}{|G|}\Bigg(2\pi n_{\rho}+\sum_{g\in G}\sqrt{3-\chi_{\rho_{1}}(g)}\overline{\chi}_{\rho}(g)x_{\lbr g\rbr}\Bigg),
\eea
where the first sum is over the non-trivial irreducible representations of $G$.
\end{Conjecture}
We have
\begin{theorem}[\cite{CCIT2}]\label{31}
Conjecture \ref{27} and \ref{28} hold for the polyhedral and binary polyhedral groups of type $A$.
\end{theorem}

\section{Induction from normal subgroups}
Let $H$ be a normal subgroup of $G$, and $i: H \hookrightarrow G$ be the injection. $i$ induces a map from the set of the congjugation classes of $H$ to the conjugation classes of $G$, which we denote by $i^{\#}: \textrm{Conj}(H) \rightarrow \textrm{Conj}(G)$. The representation $\rho$ restricts to be a representation of $H$ (for which we denote still by $\rho$), and we can consider the equivariant Gromov-Witten invariants of the corresponding $[\mathbb{C}^n /H]$. There is a simple relation between the genus 0 Gromov-Witten invariants of $[\mathbb{C}^n /G]$ with monodromies lying in $H$, and the  genus 0 Gromov-Witten invariants of $[\mathbb{C}^n /H]$. This fact was used in \cite{BGh0708}, and was also mentioned in \cite{BGh0707}. We state it as
\begin{proposition}\label{105}
\bea
F_{0}^{[\mathbb{C}^n /G]}\Big|_{t_{l}^{\lbr \gamma\rbr}=0,\gamma\not\in  Image(i^{\#}) }=\frac{|H|}{|G|}F_{0}^{[\mathbb{C}^n /H]}\Big|_{t_{l}^{\lbr \gamma\rbr}=t_{l}^{ i^{\#}(\gamma)}}.
\eea
\end{proposition}
Proof: The factor $\frac{|H|}{|G|}$ is due to the degree of $\overline{\mathcal{M}}_{g,n}(\mathcal{B}H)$ to $\overline{\mathcal{M}}_{g,n}(\mathcal{B}G)$. As in the proof of the lemma 7 in \cite{BGh0708}, we only need to show that the structure group of the corresponding $G$-torsor is reduced to $H$, when the monodromy data at every marked point is dictated in $H$. Thus it suffices to show that the generator of the local group of every nodal point lies in $H$. For smooth orbicurves this holds tautologically. It holds for nodal orbicurves because, firstly we are considering orbicurves  of genus 0, so we can do induction from the leaves of the dual \emph{tree} graph, and secondly we have the condition that the orbicurve is \emph{balanced} at every nodal point.\\

We need also to show that the Hurwitz-Hodge bundle associated to the representation $\rho$ on $\overline{\mathcal{M}}_{g,n}(\mathcal{B}G)$ pulls back to be the corresponding Hurwitz-Hodge bundle on $\overline{\mathcal{M}}_{g,n}(\mathcal{B}H)$ . Let $C$ be a balanced stable map from a marked orbicurve to $BG$.  The monodromy data at the marked points all lie in $H$ just means that the morphism $C\rightarrow BG$ factors through the canonical morphism $BH\rightarrow BG$. Thus we have the following cartesian graph
$$\xymatrix{C_{G}\ar[r]\ar[d]&\coprod\limits_{|G/H|}pt\ar[r]\ar[d] & pt\ar[d]\\
C\ar[r]& BH\ar[r]&BG    .      }
$$
From the right square one sees that the ordinary curve $C_{G}$ is a disjoint union of $|G/H|$ copies of $C_{H}:= C\underset{BH}{\times}pt$, and the quotient group $G/H$ acts on the set of copies freely and transitively. Therefore we have natural isomorphisms of cohomology groups $$H^{i}(C_{G}, \mathcal{O}_{C_{G}}\bigotimes V_{\rho}^{\vee})^{G}=\Big(\coprod\limits_{|G/H|}H^{i}(C_{H}, \mathcal{O}_{C_{H}}\bigotimes V_{\rho}^{\vee})^{H}\Big)^{G/H}=H^{i}(C_{H}, \mathcal{O}_{C_{H}}\bigotimes V_{\rho}^{\vee})^{H}.$$
\hfill\qedsymbol\\

This proposition imposes a \emph{compatibility condition} on the conjectural formulae (\ref{25}) and (\ref{26}) for genus 0 primary Gromov-Witten invariants, and we can verify the compatibility to obtain some Hurwitz-Hodge integrals inductively. We denote the right-handsides of (\ref{25}) and (\ref{26}) by $\tilde{F}_{0}^{[\mathbb{C}^2/\hat{D}_{n}]}$ and $\tilde{F}_{0}^{[\mathbb{C}^3/D_{n}]}$, respectively. We are going to verify that they satisfy \ref{105}. We shall make use of the following lemma frequently.\\

\begin{lemma}\label{29}
\ben
\sum_{k=0}^{n-1}h(x+\frac{2k\pi}{n})=\left\{\begin{array}{ll}
\frac{1}{n^2}h(nx+\pi), & 2|n , \vspace{0.4cm}, \\
\frac{1}{n^2}h(nx), & 2\nmid n.
\end{array}\right.
\een
\end{lemma}
Proof: Taking derivatives for three times, we are left to prove
\ben
\sum_{k=0}^{n-1}\tan(x+\frac{k\pi}{n})=\left\{\begin{array}{ll}
n\tan(nx+\frac{\pi}{2}), & 2|n , \vspace{0.4cm}, \\
n\tan(nx), & 2\nmid n.
\end{array}\right.
\een
But
\ben
\sum_{k=0}^{n-1}\tan(x+\frac{k\pi}{n}) &=&-\frac{d}{dx} \sum_{k=0}^{n-1}\log \cos (x+\frac{k\pi}{n})= -\frac{d}{dx}\log \prod_{k=0}^{n-1} \frac{e^{i(x+\frac{k\pi}{n})}+e^{-x-\frac{k\pi}{n}}}{2}\\
&=& -\frac{d}{dx}\log\big( e^{i\frac{n-1}{2}\pi}(e^{inx}+(-1)^{n-1}e^{-inx})\big)\\
&=&\left\{\begin{array}{ll}
-\frac{d}{dx}\log\sin (x), & 2|n , \vspace{0.3cm} \\
-\frac{d}{dx}\log\cos (x), & 2\nmid n.
\end{array}\right.
\een
\hfill\qedsymbol
\subsection{The compatibility between $\hat{D}_{n}$ and $A_{2n-5}$}

For the $\hat{D}_{n}$ generated by $a, b$ and relations  $a^{n-2}=b^{2}$, $b^{4}=1$, $ba=a^{-1}b$, where $n\geq 4$,  the subgroup generated by $a$ is a normal subgroup, which is isomorphic to the cyclic group $\mathbb{Z}_{2n-4}$, or denoted by $A_{2n-5}$. \\

For $\hat{D}_{n}$, write $x_l$ for the coordinate corresponding to $e_{\lbr a^l\rbr}$, where $1\leq l\leq n-2$, and write $y$, $z$ for the
coordinate corresponding to $e_{\lbr b\rbr}$, $e_{\lbr ab\rbr}$ respectively.
we have
\ben
P_{\rho_k}=\frac{1}{4n-8}\big ( 4\pi+2\sum_{l=1}^{n-3}\sqrt{2-\omega^{l}-\omega^{-l}}(\omega^{kl}+\omega^{-kl})x_{l}+2(\omega^{k(n-2)}+\omega^{-k(n-2)})x_{n-2}\big ),
\een
for $1\leq k\leq n-3$, and
\ben
P_{\psi_1}\Big|_{y=z=0}=P_{\psi_2}\Big|_{y=z=0}=\frac{1}{4n-8}\big ( 2\pi+2\sum_{l=1}^{n-3}\sqrt{2-\omega^{l}-\omega^{-l}}x_{l}+2x_{n-2}\big ),
\een

\ben
P_{\psi_3}\Big|_{y=z=0}=P_{\psi_4}\Big|_{y=z=0}=\frac{1}{4n-8}\big ( 2\pi+2\sum_{l=1}^{n-3}\sqrt{2-\omega^{l}-\omega^{-l}}(-1)^{l}x_{l}+2(-1)^{n}x_{n-2}\big ).
\een

For $A_{2n-5}$, write $\hat{x}_l$ for the coordinate corresponding to $e_{\lbr a^l\rbr}$, where $1\leq l\leq 2n-5$, we have
\ben
P_{\sigma_k}=\frac{1}{2n-4}\big ( 2\pi+\sum_{l=1}^{2n-5}\sqrt{2-\omega^{l}-\omega^{-l}}\omega^{-kl}\hat{x}_{l}\big ),
\een
for $1\leq k \leq 2n-5$. Replace $\hat{x}_l$ by $x_l$ for $1\leq l\leq n-2$, and $\hat{x}_l$ by $x_{2n-4-l}$ for $n-3\leq l\leq 2n-5$, we obtain
\ben
P_{\sigma_k}=P_{\sigma_{2n-4-k}}= P_{\rho_k}
\een
for $1\leq k\leq n-3$, and
\ben
P_{\sigma_0}= 2P_{\psi_1}\Big|_{y=z=0}=2P_{\psi_2}\Big|_{y=z=0},
\een

\ben
P_{\sigma_{n-2}}= 2P_{\psi_3}\Big|_{y=z=0}=2P_{\psi_4}\Big|_{y=z=0}.
\een

For simplicity of notations, we let $t_k=P_{\rho_k}$, $1\leq k\leq n-3$, and $t_{n-2}=P_{\psi_3}\Big|_{y=z=0}=P_{\psi_4}\Big|_{y=z=0}$ in this subsection. We introduce the notation $t_{i\rightarrow j}$ to stand for $\sum_{k=i}^{j}t_{k}$, when $1\leq i<j\leq n-2$. Thus
\ben
P_{\psi_1}\Big|_{y=z=0}=\pi-\sum_{k=1}^{n-2}t_{k}.
\een

Now we can write down the contribution of  the positive roots to $\frac{1}{2\lambda}F_{0}$, using the notations of  \cite{Bourbaki}, Plate I, V.
The contribution of positive roots of the form $\epsilon_{1}-\epsilon_{j}$ ($2\leq j\leq n$) is
\ben
&&h(\pi+(\pi-t_{1\rightarrow n-2}))+h(\pi+(\pi-t_{2\rightarrow n-2}))+\cdots+h(\pi+(\pi-t_{n-2}))+h(2\pi)\\
&=& \sum_{i=1}^{n-2}h(t_{i\rightarrow n-2})+h(2\pi).
\een
The contribution of positive roots of the form $\epsilon_{i}-\epsilon_{j}$ ($2\leq i< j\leq n$) is
\ben
\sum_{1\leq i\leq j\leq n-2}h(\pi+t_{i\rightarrow j}).
\een

The contribution of positive roots of the form $\epsilon_{i}+\epsilon_{n}$ ($1\leq i<  n$) is
\ben
h(2\pi)+\sum_{1\leq i\leq n-2}h(\pi+t_{i\rightarrow n-2}).
\een

The contribution of positive roots of the form $\epsilon_{1}+\epsilon_{j}$ ($2\leq j\leq n-1$) is
\ben
&&\sum_{1\leq i\leq n-2}h(\pi+(\pi+t_{i\rightarrow n-2}))\\
&=&\sum_{1\leq i\leq n-2}h(t_{i\rightarrow n-2}).
\een
The contribution of positive roots of the form $\epsilon_{i}+\epsilon_{j}$ ($2\leq i< j\leq n-1$) is
\ben
\sum_{1\leq i<j\leq n-2}h(\pi+t_{i\rightarrow j-1}+2t_{j\rightarrow n-2}).
\een
The case of positive roots of type $A_{2n-5}$ is similar. The verification of the compatibility is by using the identity $h(x)+h(\pi+x)=\frac{1}{4}h(\pi+2x)$ many times. Note that this is reminiscent of the procedure of two-step resolutions of the singularities of type D. It is also related to the partial crepant resolution conjecture \cite{CaTo}.

\subsection{The compatibility between $\hat{D}_{n}$ and $\hat{D}_{2n-2}$ }
For the $\hat{D}_{2n-2}$ generated by $a, b$ and relations  $a^{2n-4}=b^{2}$, $b^{4}=1$, $ba=a^{-1}b$, where $n\geq 4$, the subgroup generated by $a^2$ and $b$ is a normal subgroup, which is isomorphic to $\hat{D}_{n}$.\\

$V=\rho_{1}$, $\omega=\exp (\frac{2\pi \sqrt{-1}}{2n-4})$. When $n$ is even, we have
\ben
P_{\psi_1}=\frac{1}{4n-8}\big ( 2\pi+2\sum_{l=1}^{n-3}\sqrt{2-\omega^{l}-\omega^{-l}}x_{l}+2x_{n-2}+(n-2)\sqrt{2}y+(n-2)\sqrt{2}z\big ),
\een

\ben
P_{\psi_2}=\frac{1}{4n-8}\big ( 2\pi+2\sum_{l=1}^{n-3}\sqrt{2-\omega^{l}-\omega^{-l}}x_{l}+2x_{n-2}-(n-2)\sqrt{2}y-(n-2)\sqrt{2}z\big ),
\een

\ben
P_{\psi_3}=\frac{1}{4n-8}\big ( 2\pi+2\sum_{l=1}^{n-3}\sqrt{2-\omega^{l}-\omega^{-l}}(-1)^{l}x_{l}+2x_{n-2}-(n-2)\sqrt{2}y+(n-2)\sqrt{2}z\big ),
\een

\ben
P_{\psi_4}=\frac{1}{4n-8}\big ( 2\pi+2\sum_{l=1}^{n-3}\sqrt{2-\omega^{l}-\omega^{-l}}(-1)^{l}x_{l}+2x_{n-2}+(n-2)\sqrt{2}y-(n-2)\sqrt{2}z\big ).
\een

Restricting every ${x_i}$ to be $0$, we have
\ben
P_{\psi_1}=\frac{\pi}{2n-4}+\frac{\sqrt{2}}{4}y+\frac{\sqrt{2}}{4}z,\hspace{1cm} P_{\psi_2}=\frac{\pi}{2n-4}-\frac{\sqrt{2}}{4}y-\frac{\sqrt{2}}{4}z,
\een

\ben
P_{\psi_3}=\frac{\pi}{2n-4}-\frac{\sqrt{2}}{4}y+\frac{\sqrt{2}}{4}z,\hspace{1cm} P_{\psi_4}=\frac{\pi}{2n-4}+\frac{\sqrt{2}}{4}y-\frac{\sqrt{2}}{4}z.
\een

 When $n$ is odd, we have
\ben
P_{\psi_1}=\frac{1}{4n-8}\big ( 2\pi+2\sum_{l=1}^{n-3}\sqrt{2-\omega^{l}-\omega^{-l}}x_{l}+2x_{n-2}+(n-2)\sqrt{2} y+(n-2)\sqrt{2} z\big ),
\een

\ben
P_{\psi_2}=\frac{1}{4n-8}\big ( 2\pi+2\sum_{l=1}^{n-3}\sqrt{2-\omega^{l}-\omega^{-l}}x_{l}+2x_{n-2}-(n-2)\sqrt{2}y-(n-2)\sqrt{2}z\big ),
\een

\ben
P_{\psi_3}=\frac{1}{4n-8}\big ( 2\pi+2\sum_{l=1}^{n-3}\sqrt{2-\omega^{l}-\omega^{-l}}(-1)^{l}x_{l}-2x_{n-2}+(n-2)\sqrt{2}i y-(n-2)\sqrt{2}i z\big ),
\een

\ben
P_{\psi_4}=\frac{1}{4n-8}\big ( 2\pi+2\sum_{l=1}^{n-3}\sqrt{2-\omega^{l}-\omega^{-l}}(-1)^{l}x_{l}-2x_{n-2}-(n-2)\sqrt{2}i y+(n-2)\sqrt{2}i z\big ).
\een

Restricting every ${x_i}$ to be $0$, we have
\ben
P_{\psi_1}=\frac{\pi}{2n-4}+\frac{\sqrt{2}}{4}y+\frac{\sqrt{2}}{4}z,\hspace{1cm} P_{\psi_2}=\frac{\pi}{2n-4}-\frac{\sqrt{2}}{4}y-\frac{\sqrt{2}}{4}z,
\een

\ben
P_{\psi_3}=\frac{\pi}{2n-4}+\frac{\sqrt{2}}{4}i y-\frac{\sqrt{2}}{4}i z,\hspace{1cm} P_{\psi_4}=\frac{\pi}{2n-4}-\frac{\sqrt{2}}{4}i y+\frac{\sqrt{2}}{4}i z.
\een

When $n$ is even, we have
\ben
&&\frac{1}{2\lambda}\tilde{F}_{0}^{[\mathbb{C}^2 /\hat{D}_n]}\Big |_{\forall x_i =0}\\
&=&\sum_{k=0}^{n-3}h\Big(\frac{(2n-3)\pi}{2n-4}-\frac{\sqrt{2}}{4}(y+z)+\frac{k\pi}{n-2}\Big )+h(\frac{\sqrt{2}}{2}y)\\
&&+h\Big(\frac{\pi}{2n-4}-\frac{\sqrt{2}}{4}(y+z)\Big)+h(\frac{\sqrt{2}}{2}z)+\sum_{k=0}^{n-3}h\Big(\frac{(2n-3)\pi}{2n-4}-
\frac{\sqrt{2}}{4}(y-z)+\frac{k\pi}{n-2}\Big )\\
&&+\sum_{k=0}^{n-3}h\Big(\frac{(2n-3)\pi}{2n-4}+\frac{\sqrt{2}}{4}(y-z)+\frac{k\pi}{n-2}\Big )+\sum_{k=0}^{n-3}h\Big(\frac{(2n-3)\pi}{2n-4}-\frac{\sqrt{2}}{4}(y+z)+\frac{(2n-4-k)\pi}{n-2}\Big )\\
&=&\sum_{k=0}^{2n-5}h\Big(\frac{(2n-3)\pi}{2n-4}+\frac{\sqrt{2}}{4}(y-z)+\frac{k\pi}{n-2}\Big )+\sum_{k=0}^{2n-5}h\Big(\frac{(2n-3)\pi}{2n-4}-\frac{\sqrt{2}}{4}(y+z)+\frac{k\pi}{n-2}\Big)\\
&&+h(\frac{\sqrt{2}}{2}y)+h(\frac{\sqrt{2}}{2}z)\\
&=& \frac{1}{4(n-2)^2}h\Big(\frac{\sqrt{2}(n-2)}{2}(y-z)\Big)+\frac{1}{4(n-2)^2}h\Big(\frac{\sqrt{2}(n-2)}{2}(y+z)\Big )+h(\frac{\sqrt{2}}{2}y)+h(\frac{\sqrt{2}}{2}z).
\een
When $n$ is odd, a similar computation gives,
\ben
&&\frac{1}{2t}\tilde{F}_{0}^{[\mathbb{C}^2 /\hat{D}_n]}\Big |_{\forall x_i =0}\\
&&=\frac{1}{4(n-2)^2}h\Big(\frac{\sqrt{2}(n-2)}{2}i(y-z)\Big)+\frac{1}{4(n-2)^2}h\Big(\frac{\sqrt{2}(n-2)}{2}(y+z)\Big )\\
&&+h(\frac{\sqrt{2}}{2}y)+h(\frac{\sqrt{2}}{2}z).
\een
Thus it is easy to see that
\ben
\tilde{F}^{[\mathbb{C}^2 /\hat{D}_{2n}]}\Big |_{\forall x_i =0, z=0}=\frac{1}{2}\tilde{F}^{[\mathbb{C}^2 /\hat{D}_{n+1}]}\Big |_{\forall x_i =0, z=y}.
\een

\begin{remark}
Let $n$ be even, and take $z=0$. We have
\ben
&&\frac{1}{2\lambda}\sum_{l\geq 3}\frac{y^{l-3}}{(l-3)!}\langle e_{\lbr b \rbr}^{l} \rangle^{[\mathbb{C}^2 /\hat{D}_n]}\\
&=&\frac{d^3}{dy^3}\Big(\frac{1}{2\lambda}F^{[\mathbb{C}^2 /\hat{D}_{2n}]}\Big |_{\forall x_i =0, z=0}\Big)\\
&=& \frac{n-2}{4\sqrt{2}}\tan (-\frac{n-2}{2\sqrt{2}}y)+\frac{1}{2\sqrt{2}}\tan (-\frac{y}{2\sqrt{2}})\\
&=& \frac{n-2}{4}\sum_{m=1}^{\infty}\frac{(-1)^{m}(2^{2m}-1)B_{2m}(n-2)^{2m-1}}{2^{2m-1}(2m)!}y^{2m-1}\\
&& +\frac{1}{2}\sum_{m=1}^{\infty}\frac{(-1)^{m}B_{2m}}{2^{m-1}(2m)!}y^{2m-1}.
\een
Thus
\ben
\langle e_{\lbr b \rbr}^{2m} \rangle^{[\mathbb{C}^2 /\hat{D}_n]}=\lambda\Big(\frac{(n-2)^{2m-2}}{2}+1\Big)\frac{(-1)^{m-1}(2^{2m-2}-1)B_{2m-2}}{2^{m-1}(2m-2)!}.
\een
When $m=2$, this coincides with the result of example \ref{19}. Note that a discrepancy of sign arises because here we are integrating $(R\pi_{*}^{1}{C}_{\hat{D}_{n}}-R\pi_{*}^{0}{C}_{\hat{D}_{n}})^{\rho_1}$, by definition.
\end{remark}

The following lemma is easy, we record it here for later use.
\begin{lemma}\label{65}
\ben
\tilde{F}_{0}^{[\mathbb{C}^{2}/\hat{D}_n]}(x_0,\cdots,x_{n-2},y,z)=\tilde{F}_{0}^{[\mathbb{C}^{2}/\hat{D}_n]}(x_0,\cdots,x_{n-2},z,y).
\een
\end{lemma}
Proof: This follows from the fact that there is an automorphism of $\hat{D}_{n}$ sending $a$ to $a$, and $b$ to $ab$, while the representation $\rho_{1}$ is equivalent to itself by this automorphism. \hfill\qed

\section{Induction by the WDVV equation}

\subsection{The case of $[\mathbb{C}^3 /D_n]$ }
As a warming up and also for later use, we first show the WDVV-induction for the $D_n\subset SO(3)$. The group $D_n$ is generated by $a,b$ and the relations $a^{n-2}=b^2=1, ba=a^{-1}b$. We list  the conjugacy classes, the elements of the conjugacy classes, the cardinality of  the centralizer of an element in the conjugacy classes and the corresponding variables we use in the generating function as the following. For Let $z(g)$ be the order of the centralizer of an element $g$, and $z(\lbr g \rbr)$ the order of the centralizer of an element in the conjugacy class $\lbr g\rbr$. When $n$ is odd, the table is\\
 \begin{center}
\begin{tabular}{l|l|l|l}\hline
conjugacy classes & $\lbr 1 \rbr$ & $ \lbr a^k \rbr $, $1\leq k \leq \frac{n-3}{2} $ & $\lbr b \rbr $  \\ \hline
elements          & $\{1\}$           & $\{a^k, a^{n-2-k}\}$, $1\leq k\leq n-3 $ &$\{b, ab, \cdots, a^{n-3}b\}$  \\ \hline
$z(\lbr g\rbr)$ & $2n-4$ & $n-2$ & $2$ \\ \hline
variables         & $x_0 $        & $x_k$, $1\leq k\leq \frac{n-3}{2}$ & $y$ \\ \hline
\end{tabular}
\end{center}
We use $\langle\cdot\rangle^{G}$ to denote the genus zero Gromov-Witten  invariants of $BG$. The only nonzero length three correlators are
\ben
\langle e_{\lbr a^k \rbr} e_{\lbr b\rbr}e_{\lbr b\rbr}\rangle ^{{D}_n}=1,
\een
\ben
\langle e_{\lbr a^{i}\rbr} e_{\lbr a^{j} \rbr}e_{\lbr a^{k}\rbr}\rangle ^{{D}_n}=\frac{1}{n-2}, \hspace{0.2cm}\text{for}\hspace{0.2cm}1\leq i,j,k\leq \frac{n-3}{2} \hspace{0.2cm}\text{satisfying}\hspace{0.2cm} \pm i \pm j\pm k\equiv 0 (\text{mod}\hspace{0.2cm}n-2),
\een
\ben
\langle e_{\lbr 1\rbr} e_{\lbr a^k \rbr} e_{\lbr a^k\rbr}\rangle ^{{D}_n}=\frac{1}{n-2},  \quad 1\leq k\leq n-3,
\een
\ben
\langle e_{\lbr 1\rbr} e_{\lbr b \rbr} e_{\lbr b\rbr}\rangle ^{{D}_n}=\frac{1}{2},
\een
\ben
\langle e_{\lbr 1\rbr} e_{\lbr 1 \rbr} e_{\lbr 1\rbr}\rangle ^{{D}_n}=\frac{1}{2n-4}.
\een
Note that for the dimensional reason, the length \emph{three} correlators $\langle\cdot\rangle^{[\mathbb{C}^3 /D_n]}$ without $e_{\lbr 1\rbr}$ insertions are related to those of $\mathcal{B}{D}_n$ by
\bea\label{61}
\langle\cdot\rangle^{[\mathbb{C}^3 /D_n]}=\langle\cdot\rangle^{D_n},
\eea
also, for length three correlator with exactly one $e_{\lbr 1\rbr}$ insertion we have
\bea\label{62}
\langle\cdot\rangle^{[\mathbb{C}^3 /D_n]}=\lambda^{-1}\langle\cdot\rangle^{D_n}.
\eea

This and the next subsection are based on the following
\begin{proposition}\label{63}
The genus 0 primary equivariant Gromov-Witten potential $F_{0}^{[\mathbb{C}^2 /\hat{D}_n]}$ and $F_{0}^{[\mathbb{C}^3 /{D}_n]}$ satisfies the axiom of fundamental class (the string equation) and the WDVV equations.
\end{proposition}
Proof: This follows from the theorem \ref{18}, and Givental's \emph{symplectic geometry of Frobenius structures} \cite{Gi}. Note that the coefficient $g^{ij}$ in the WDVV equations is diagonal in the cases we concern.
\hfill\qedsymbol
\\

In the following part in this subsection, we write $\langle\cdot\rangle$ for $\langle\cdot\rangle^{[\mathbb{C}^3 /D_n]}$ for simplicity of notations. We follow the strategy of \cite{BGh0708}, i.e., we select four insertions $e_{\lbr \gamma_{i}\rbr}$, $e_{\lbr \gamma_{j}\rbr}$, $e_{\lbr \gamma_{k}\rbr}$, $e_{\lbr \gamma_{l}\rbr}$, and expand the corresponding WDVV equation, thus obtain an equality of two generating functions. The equality of coefficients will given many equations of correlators. For one of the equations, we treat the involved correlators of the maximal length as unknown number, and those of less length as known, thus we obtain a \emph{linear} equations of some correlators of the same length. If we can find sufficiently many independent such linear equations, we can determine some correlators from more simple ones.
Following \cite{BGh0708}, we use the notation $\langle \cdots(e_{\lbr \gamma_i\rbr}e_{\lbr \gamma_j\rbr}|e_{\lbr \gamma_k\rbr}e_{\lbr \gamma_l\rbr})\rangle=\langle \cdots(e_{\lbr \gamma_i\rbr}e_{\lbr \gamma_k\rbr}|e_{\lbr \gamma_j\rbr}e_{\lbr \gamma_l\rbr})\rangle$ to indicate that we are expanding the WDVV \bea\label{106}
F_{e_{\lbr \gamma_i\rbr}e_{\lbr \gamma_j\rbr}e_{\lbr\alpha\rbr}}g^{e_{\lbr\alpha\rbr}e_{\lbr\beta\rbr}}F_{e_{\lbr\beta\rbr} e_{\lbr \gamma_k\rbr}e_{\lbr \gamma_l\rbr}}=F_{e_{\lbr \gamma_i\rbr}e_{\lbr \gamma_k\rbr}e_{\lbr\alpha\rbr}}g^{e_{\lbr\alpha\rbr}e_{\lbr\beta\rbr}}F_{e_{\lbr\beta\rbr} e_{\lbr \gamma_j\rbr}e_{\lbr \gamma_l\rbr}}.
\eea
For example, consider the WDVV equation of the form
\bea\label{101}
&&\langle e_{\lbr a\rbr }^{k_1}\cdots  e_{\lbr a^{\frac{n-3}{2}}\rbr }^{k_{\frac{n-3}{2}}} e_{\lbr b\rbr }^{k_{b}} ( e_{\lbr a^{i}\rbr } e_{\lbr b\rbr }| e_{\lbr a^{i}\rbr } e_{\lbr b\rbr })\rangle\nonumber\\
&=&\langle e_{\lbr a\rbr }^{k_1}\cdots  e_{\lbr a^{\frac{n-3}{2}}\rbr }^{k_{\frac{n-3}{2}}} e_{\lbr b\rbr }^{k_{b}} ( e_{\lbr a^{i}\rbr } e_{\lbr a^{i}\rbr }| e_{\lbr b\rbr } e_{\lbr b\rbr })\rangle,
\eea
where $k_1+\cdots +k_{\frac{n-3}{2}}+k_{b}=m-3\geq 1$, $1\leq i\leq \frac{n-3}{2}$. Expand the left hand side, we obtain
\bea\label{101}
&&\langle e_{\lbr a\rbr }^{k_1}\cdots  e_{\lbr a^{i}\rbr }^{k_i+1}\cdots  e_{\lbr a^{\frac{n-3}{2}}\rbr }^{k_{\frac{n-3}{2}}} e_{\lbr b\rbr }^{k_{b}+2}\rangle \big(\langle e_{\lbr 1\rbr} e_{\lbr b \rbr} e_{\lbr b\rbr}\rangle \big)^{-1}\langle e_{\lbr a^{i}\rbr } e_{\lbr b\rbr }  e_{\lbr b\rbr }\rangle\nonumber\\
&&+\langle e_{\lbr a^{i}\rbr } e_{\lbr b\rbr }  e_{\lbr b\rbr }\rangle  \big(\langle e_{\lbr 1\rbr} e_{\lbr b \rbr} e_{\lbr b\rbr}\rangle \big)^{-1}\langle e_{\lbr a\rbr }^{k_1}\cdots  e_{\lbr a^{i}\rbr }^{k_i+1}\cdots  e_{\lbr a^{\frac{n-3}{2}}\rbr }^{k_{\frac{n-3}{2}}} e_{\lbr b\rbr }^{k_{b}+2}\rangle+\text{Length}(<m)\nonumber \\
&=& 4\lambda\langle e_{\lbr a\rbr }^{k_1}\cdots  e_{\lbr a^{i}\rbr }^{k_i+1}\cdots  e_{\lbr a^{\frac{n-3}{2}}\rbr }^{k_{\frac{n-3}{2}}} e_{\lbr b\rbr }^{k_{b}+2}\rangle+\text{Length}(<m).
\eea
Next we expand the right hand side. To simplify the expressions, we introduce a map $\phi:\{1,\cdots,\frac{n-3}{2}\}\rightarrow\{1,\cdots,\frac{n-3}{2}\}$, such that
$$
\phi(i)=\left\{\begin{array}{ll}
2i, & \text{if}\quad 2i\leq \frac{n-3}{2},\\
n-2-2i, & \text{if}\quad  2i> \frac{n-3}{2}.
\end{array}\right.
$$
Since $n$ is odd, $\phi$ is \emph{bijective}. Thus the right hand side is equal to
\bea\label{102}
&&\sum_{1\leq l\leq \frac{n-3}{2}}\langle e_{\lbr a\rbr }^{k_1}\cdots  e_{\lbr a^{i}\rbr }^{k_i+2}\cdots e_{\lbr a^{l}\rbr }^{k_{l}+1}\cdots  e_{\lbr a^{\frac{n-3}{2}}\rbr }^{k_{\frac{n-3}{2}}} e_{\lbr b\rbr }^{k_{b}}\rangle \big(\langle e_{\lbr 1\rbr} e_{\lbr a^{l} \rbr} e_{\lbr a^{l}\rbr}\rangle \big)^{-1}\langle e_{\lbr a^{l}\rbr } e_{\lbr b\rbr }  e_{\lbr b\rbr }\rangle\nonumber\\
&&+\langle e_{\lbr a^{i}\rbr } e_{\lbr a^{i}\rbr }  e_{\lbr a^{\phi(i)}\rbr }\rangle  \big(\langle e_{\lbr 1\rbr} e_{\lbr a^{\phi(i)} \rbr} e_{\lbr a^{\phi(i)}\rbr}\rangle \big)^{-1}\langle e_{\lbr a\rbr }^{k_1}\cdots  e_{\lbr a^{\phi(i)}\rbr }^{k_{\phi(i)}+1}\cdots  e_{\lbr a^{\frac{n-3}{2}}\rbr }^{k_{\frac{n-3}{2}}} e_{\lbr b\rbr }^{k_{b}+2}\rangle+\text{Length}(<m)\nonumber \\
&=& (n-2)\lambda\sum_{1\leq l\leq \frac{n-3}{2}}\langle e_{\lbr a\rbr }^{k_1}\cdots  e_{\lbr a^{i}\rbr }^{k_i+2}\cdots e_{\lbr a^{l}\rbr }^{k_{l}+1}\cdots  e_{\lbr a^{\frac{n-3}{2}}\rbr }^{k_{\frac{n-3}{2}}} e_{\lbr b\rbr }^{k_{b}}\rangle \nonumber\\
&&+ \lambda\langle e_{\lbr a\rbr }^{k_1}\cdots  e_{\lbr a^{\phi(i)}\rbr }^{k_{\phi(i)}+1}\cdots  e_{\lbr a^{\frac{n-3}{2}}\rbr }^{k_{\frac{n-3}{2}}} e_{\lbr b\rbr }^{k_{b}+2}\rangle+\text{Length}(<m).
\eea
Set $u_{i}=\lambda\langle e_{\lbr a\rbr }^{k_1}\cdots  e_{\lbr a^{i}\rbr }^{k_i+1}\cdots  e_{\lbr a^{\frac{n-3}{2}}\rbr }^{k_{\frac{n-3}{2}}} e_{\lbr b\rbr }^{k_{b}+2}\rangle$ for $1\leq i\leq \frac{n-3}{2}$ temporarily. Arranging the equality (\ref{101})=(\ref{102}) we obtain
\bea{\label{103}}
&&4u_{i}-u_{\phi(i)}=(n-2)\lambda\sum_{1\leq l\leq \frac{n-3}{2}}\langle e_{\lbr a\rbr }^{k_1}\cdots  e_{\lbr a^{i}\rbr }^{k_i+2}\cdots e_{\lbr a^{l}\rbr }^{k_{l}+1}\cdots  e_{\lbr a^{\frac{n-3}{2}}\rbr }^{k_{\frac{n-3}{2}}} e_{\lbr b\rbr }^{k_{b}}\rangle +\text{Length}(<m).\nonumber\\
&&
\eea
  We view this as a  linear system for the variables $u_{i}$, $1\leq i\leq \frac{n-3}{2}$. Because $\phi$ is a bijection, the set $\{1,\cdots,\frac{n-3}{2}\}$ decomposes into several subsets, each subset being an orbit of the map $\phi$. Let $\{i_1,\cdots,i_l\}$ be such an orbit, thus it forms a \emph{cycle} under the iterations of $\phi$. For the corresponding variables $u_{i_1},\cdots,u_{i_l}$, the equations in (\ref{103}) involving them form a linear subsystem. The corresponding  matrix of coefficients is of the form
  $$
\left (\begin{array}{ccccc}
4& -1 & & &\\
& 4 & -1 & &  \\
&& \cdots & \cdots & \\
&&& 4 & -1 \\
-1 &&&& 4 \\
\end{array}
\right),
$$
which is easily seen to be nonsingular. Thus in this way we  can determine a correlator with at least two $e_{\lbr b\rbr}$-insertions and at least one insertions of the form $e_{\lbr a^{i}\rbr}$ from correlators of less lengths and correlators of less $e_{\lbr b\rbr}$-insertions by WDVV equations. Note that a nonzero correlator has an even number of insertions $e_{\lbr b\rbr}$, thus at least two, if there is any. Therefore we finally come to
\begin{theorem}\label{104}
For odd $n\geq 5$, The correlators $\langle \cdot \rangle ^{[\mathbb{C}^3/D_{n}]}$ of length at least four are determined by correlators of length three, the correlators with only insertions of the form $e_{\lbr a^k \rbr}$, $1\leq k\leq \frac{n-3}{2}$, the correlators $\langle e_{\lbr b \rbr}^k \rangle$ and the WDVV equations.\hfill\qedsymbol
\end{theorem}
When $n$ is even, $b$ and $ab$ are not in the same conjugate classes. But the similar theorem still holds.
\begin{theorem}\label{104}
For even $n\geq 6$, The correlators $\langle \cdot \rangle ^{[\mathbb{C}^3/D_{n}]}$ of length at least four are determined by correlators of length three, the correlators with only insertions of the form $e_{\lbr a^k \rbr}$, $1\leq k\leq \frac{n-3}{2}$, the correlators $\langle e_{\lbr b \rbr}^k \rangle$, $\langle e_{\lbr ab \rbr}^k \rangle$  and the WDVV equations.
\end{theorem}
Proof: The proof is more complicated than the former one but less complicated than and very similar to the proofs in the next subsection, so we omit the details here.\hfill\qedsymbol

\subsection{The case of $[\mathbb{C}^2 /\hat{D}_n]$}
We have the table
 \begin{center}
\begin{tabular}{l|l|l|l}\hline
conjugacy classes & $\lbr 1 \rbr$ & $ \lbr a^k \rbr $, $1\leq k \leq n-3 $ & $\lbr a^{n-2} \rbr $  \\ \hline
elements          & $\{1\}$           & $\{a^k, a^{2n-4-k}\}$, $1\leq k\leq n-3 $ &$\{a^{n-2}\}$  \\ \hline
$z(\lbr g\rbr)$ & $4n-8$ & $2n-4$ & $4n-8$ \\ \hline
variables         & $x_0 $        & $x_k$, $1\leq k\leq n-3$ & $x_{n-2}$ \\ \hline
\end{tabular}
\end{center}

\begin{center}\begin{tabular}{l|l|l}\hline
conjugacy classes &  $\lbr b \rbr $ & $\lbr ab \rbr$ \\ \hline
elements          & $\{b, a^2 b, \cdots, a^{2n-6}b\}$& $\{ab, a^3 b, \cdots, a^{2n-5}b\}$\\ \hline
$z(\lbr g\rbr)$ &  $4$ & $4$\\ \hline
variables         &  $y$ & $z$\\ \hline
\end{tabular}
\end{center}

The only nonzero length three correlators are
\ben
\langle e_{\lbr a^k \rbr} e_{\lbr b\rbr}e_{\lbr b\rbr}\rangle ^{\hat{D}_n}=\langle e_{\lbr a^k \rbr} e_{\lbr ab\rbr}e_{\lbr ab\rbr}\rangle ^{\hat{D}_n}=\frac{1}{2},\quad \text{for} \quad 1\leq k\leq n-3, \quad2|k,
\een

\ben
\langle e_{\lbr a^k \rbr} e_{\lbr b\rbr}e_{\lbr ab\rbr}\rangle ^{\hat{D}_n}=\frac{1}{2},\quad \text{for} \quad 1\leq k\leq n-3, \quad2\nmid k,
\een

\ben
\langle e_{\lbr a^{n-2} \rbr} e_{\lbr b\rbr}e_{\lbr b\rbr}\rangle ^{\hat{D}_n}=\langle e_{\lbr a^{n-2} \rbr} e_{\lbr ab\rbr}e_{\lbr ab\rbr}\rangle ^{\hat{D}_n}=\frac{1}{4},\quad \text{when} \quad 2|n,
\een

\ben
\langle e_{\lbr a^{n-2} \rbr} e_{\lbr b\rbr}e_{\lbr ab\rbr}\rangle ^{\hat{D}_n}=\frac{1}{4},\quad \text{when}  \quad2\nmid n,
\een

\ben
\langle e_{\lbr 1 \rbr} e_{\lbr b\rbr}e_{\lbr b\rbr}\rangle ^{\hat{D}_n}=\langle e_{\lbr 1 \rbr} e_{\lbr ab\rbr}e_{\lbr ab\rbr}\rangle ^{\hat{D}_n}=\frac{1}{4},
\een

\ben
\langle e_{\lbr 1\rbr} e_{\lbr 1 \rbr}e_{\lbr 1\rbr}\rangle ^{\hat{D}_n}=\langle e_{\lbr 1\rbr} e_{\lbr a^{n-2} \rbr}e_{\lbr a^{n-2}\rbr}\rangle ^{\hat{D}_n}=\frac{1}{4n-8},
\een

\ben
\langle e_{\lbr 1\rbr} e_{\lbr a^{k} \rbr}e_{\lbr a^{k}\rbr}\rangle ^{\hat{D}_n}=\frac{1}{2n-4}, \quad 1\leq k\leq n-3,
\een
and
\ben
\langle e_{\lbr a^{i}\rbr} e_{\lbr a^{j} \rbr}e_{\lbr a^{k}\rbr}\rangle ^{\hat{D}_n}=\frac{1}{2n-4}, \hspace{0.2cm}\text{for}\hspace{0.2cm}1\leq i,j,k\leq n-2 \hspace{0.2cm}\text{satisfying}\hspace{0.2cm} \pm i \pm j\pm k\equiv 0 (\text{mod}\hspace{0.2cm}2n-4).
\een

Note that for the dimensional reason, the length \emph{three} correlators $\langle\cdot\rangle^{[\mathbb{C}^2 /\hat{D}_n]}$ without $e_{\lbr 1\rbr}$ insertions are related to those of $\mathcal{B}\hat{D}_n$ by
\bea\label{61}
\langle\cdot\rangle^{[\mathbb{C}^2 /\hat{D}_n]}=\lambda\langle\cdot\rangle^{\hat{D}_n},
\eea
also, for length three correlator with exactly one $e_{\lbr 1\rbr}$ insertion we have
\bea\label{62}
\langle\cdot\rangle^{[\mathbb{C}^2 /\hat{D}_n]}=\langle\cdot\rangle^{\hat{D}_n}.
\eea

\begin{lemma}\label{1}
 The length $m(\geq 4)$ correlators with at least one insertion of the form $e_{\lbr a^k \rbr}$, $1\leq k\leq n-2$, are uniquely determined by the correlators of length less than $m$ and the correlators with only insertions of the form $e_{\lbr a^k \rbr}$, $1\leq k\leq n-2$ and the WDVV equations.
\end{lemma}
Proof: By the string equation we can only consider the correlators without insertions of $e_{\lbr 1\rbr}$. We write $\langle \cdot\rangle$ for
$\langle \cdot\rangle^{[\mathbb{C}^2 /\hat{D}_n]}$, for simplicity of notations.

For the monodromy reason, $\langle e_{\lbr a\rbr }^{k_1}\cdots  e_{\lbr a^{n-2}\rbr }^{k_{n-2}} e_{\lbr b\rbr }^{k_{n-1}} e_{\lbr ab\rbr }^{k_{n}}\rangle\neq 0$ forces $k_{n-1}+k_{n}$ to be even. We prove the lemma by induction on $k_{n-1}+k_{n}$. Thus we assume $k_{n-1}+k_{n}=2k+2$, where $k\geq 0$.\\
  Consider the WDVV equation
\ben
\langle e_{\lbr a\rbr }^{k_1}\cdots  e_{\lbr a^{n-2}\rbr }^{k_{n-2}} e_{\lbr b\rbr }^{k_{n-1}} e_{\lbr ab\rbr }^{k_{n}}( e_{\lbr a^{i}\rbr } e_{\lbr b\rbr }| e_{\lbr ab\rbr } e_{\lbr ab\rbr })\rangle=\langle e_{\lbr a\rbr }^{k_1}\cdots  e_{\lbr a^{n-2}\rbr }^{k_{n-2}} e_{\lbr b\rbr }^{k_{n-1}} e_{\lbr ab\rbr }^{k_{n}}( e_{\lbr a^{i}\rbr } e_{\lbr ab\rbr }| e_{\lbr b\rbr } e_{\lbr ab\rbr })\rangle,
\een
for any odd $i$. Expanding both sides, we obtain
\bea\label{11}
&&\sum_{1\leq l\leq n-2, 2|l}\langle \cdots e_{\lbr a^i \rbr }^{k_{i} +1}\cdots  e_{\lbr a^{l}\rbr }^{k_{l}+1} \cdots e_{\lbr b\rbr }^{k_{n-1}+1} e_{\lbr ab\rbr }^{k_{n}}\rangle\cdot (n-2)\nonumber\\
&&+2\langle \cdots e_{\lbr b\rbr }^{k_{n-1}} e_{\lbr ab\rbr }^{k_{n}+3}\rangle+\text{Length}(<m)\nonumber \\
&=&\sum_{1\leq l\leq n-2, 2\nmid l}\langle \cdots e_{\lbr a^i \rbr }^{k_{i} +1}\cdots  e_{\lbr a^{l}\rbr }^{k_{l}+1} \cdots e_{\lbr b\rbr }^{k_{n-1}} e_{\lbr ab\rbr }^{k_{n}+1}\rangle\cdot (n-2)\nonumber\\
&&+2\langle \cdots e_{\lbr b\rbr }^{k_{n-1}+2} e_{\lbr ab\rbr }^{k_{n}+1}\rangle+\text{Length}(<m),
\eea
here we follow the notations in \cite{BGh0708}, to let $\text{Length}(<m)$ stand for any combinations (allowing products and sums) of correlators of length less than $m$. Alternating the roles of $\langle e_{\lbr b \rbr} \rangle$ and $\langle e_{\lbr ab \rbr} \rangle$ we have
\bea\label{12}
&&\sum_{1\leq l\leq n-2, 2|l}\langle \cdots e_{\lbr a^i \rbr }^{k_{i} +1}\cdots  e_{\lbr a^{l}\rbr }^{k_{l}+1} \cdots e_{\lbr b\rbr }^{k_{n-1}} e_{\lbr ab\rbr }^{k_{n}+1}\rangle\cdot (n-2)\nonumber\\
&&+2\langle \cdots e_{\lbr b\rbr }^{k_{n-1}+3} e_{\lbr ab\rbr }^{k_{n}}\rangle+\text{Length}(<m) \nonumber\\
&=&\sum_{1\leq l\leq n-2, 2\nmid l}\langle \cdots e_{\lbr a^i \rbr }^{k_{i} +1}\cdots  e_{\lbr a^{l}\rbr }^{k_{l}+1} \cdots e_{\lbr b\rbr }^{k_{n-1}+1} e_{\lbr ab\rbr }^{k_{n}}\rangle\cdot (n-2)\nonumber\\
&&+2\langle \cdots e_{\lbr b\rbr }^{k_{n-1}+1} e_{\lbr ab\rbr }^{k_{n}+2}\rangle+\text{Length}(<m).
\eea
From (\ref{11}) and (\ref{12}) and lemma \ref{65} we see that, when $k_{n-1}+k_{n}=2k+2$ is fixed, we need only consider the correlators of the form
\ben
\langle\cdots e_{\lbr b\rbr }^{k+1} e_{\lbr ab\rbr }^{k+1}\rangle
\een
or
\ben
\langle\cdots e_{\lbr b\rbr }^{k+2} e_{\lbr ab\rbr }^{k}\rangle.
\een
From now on we treat the two cases $2|n$ and $2\nmid n$ separately for clearness. First we assume $n$ is even.\\
Consider the WDVV equation of the form
\bea\label{2}
&&\langle e_{\lbr a\rbr }^{k_1}\cdots  e_{\lbr a^{n-2}\rbr }^{k_{n-2}} e_{\lbr b\rbr }^{k_{n-1}} e_{\lbr ab\rbr }^{k_{n}}( e_{\lbr a^{i}\rbr } e_{\lbr a^{i}\rbr }| e_{\lbr b\rbr } e_{\lbr ab\rbr })\rangle\nonumber\\
&=&\langle e_{\lbr a\rbr }^{k_1}\cdots  e_{\lbr a^{n-2}\rbr }^{k_{n-2}} e_{\lbr b\rbr }^{k_{n-1}} e_{\lbr ab\rbr }^{k_{n}}( e_{\lbr a^{i}\rbr } e_{\lbr b\rbr }| e_{\lbr a^{i}\rbr } e_{\lbr ab\rbr })\rangle,
\eea
where $k_1+\cdots +k_{n}=m-3\geq 1$, $1\leq i\leq n-3$. Expand the left hand side, we obtain
\ben
\text{LHS}&=&\sum_{1\leq l\leq n-2}\langle e_{\lbr a\rbr }^{k_1}\cdots e_{\lbr a^{i-1}\rbr }^{k_{i-1}} e_{\lbr a^{i}\rbr }^{k_{i}+2}  e_{\lbr a^{i+1}\rbr }^{k_{i+1}}\cdots e_{\lbr a^{l-1}\rbr }^{k_{l-1}}e_{\lbr a^{l}\rbr }^{k_{l}+1}e_{\lbr a^{l+1}\rbr }^{k_{l+1}} \cdots e_{\lbr b\rbr }^{k_{n-1}} e_{\lbr ab\rbr }^{k_{n}}\rangle \\
&&z(\lbr a^l \rbr)\langle  e_{\lbr a^{l}\rbr }e_{\lbr b\rbr }e_{\lbr ab\rbr }\rangle +\sum_{1\leq l\leq n-2}\langle e_{\lbr a\rbr }^{k_1}\cdots e_{\lbr a^{l-1}\rbr }^{k_{l-1}} e_{\lbr a^{l}\rbr }^{k_{l}+1}  e_{\lbr a^{l+1}\rbr }^{k_{l+1}}\cdots e_{\lbr b\rbr }^{k_{n-1}+1} e_{\lbr ab\rbr }^{k_{n}+1}\rangle \\
&&z(\lbr a^l \rbr)\langle  e_{\lbr a^{l}\rbr }e_{\lbr a^i\rbr }e_{\lbr a^i\rbr }\rangle+\text{Length}(<m),
\een
 Now we use the proceeding computation of correlators of length $3$. To simplify the expressions, we introduce a map $\phi:\{1,\cdots,n-3\}\rightarrow\{1,\cdots,n-3,n-2\}$, such that
$$
\phi(i)=\left\{\begin{array}{ll}
2i, & \text{if}\quad 1\leq i\leq \frac{n-2}{2},\\
2n-4-2i, & \text{if}\quad \frac{n-2}{2}< i\leq n-3.
\end{array}\right.
$$
 Thus we have
\bea\label{7}
\text{LHS}&=&\sum_{1\leq l\leq n-2, 2\nmid l}\langle\cdots e_{\lbr a^{i}\rbr }^{k_{i}+2}\cdots e_{\lbr a^{l}\rbr }^{k_{l}+1}\cdots \rangle \cdot (n-2)\lambda\nonumber\\
&&+\langle \cdots  e_{\lbr a^{\phi(i)}\rbr }^{k_{\phi(i)}+1} \cdots e_{\lbr b\rbr }^{k_{n-1}+1} e_{\lbr ab\rbr }^{k_{n}+1}\rangle\cdot\big(1\cdot(1-\delta_{i,\frac{n-2}{2}})\nonumber\\
&&+2\delta_{i,\frac{n-2}{2}}\big)\lambda+\text{Length}(<m).
\eea
Expand the right hand side of (\ref{2}) as
\bea\label{8}
\text{RHS}&=&\langle \cdots e_{\lbr a^{i}\rbr }^{k_{i}+1}   \cdots e_{\lbr b\rbr }^{k_{n-1}+2} e_{\lbr ab\rbr }^{k_{n}}\rangle z(\lbr b \rbr)\langle  e_{\lbr b\rbr }e_{\lbr a^i \rbr }e_{\lbr ab\rbr }\rangle +\langle \cdots e_{\lbr a^{i}\rbr }^{k_{i}+1}   \cdots e_{\lbr b\rbr }^{k_{n-1}+1} e_{\lbr ab\rbr }^{k_{n}+1}\rangle\nonumber\\
&& z(\lbr ab \rbr)\langle  e_{\lbr ab\rbr }e_{\lbr a^i \rbr }e_{\lbr ab\rbr }\rangle +\langle  e_{\lbr b\rbr }e_{\lbr a^i \rbr }e_{\lbr ab\rbr }\rangle z(\lbr ab\rbr)\langle \cdots e_{\lbr a^{i}\rbr }^{k_{i}+1}   \cdots e_{\lbr ab\rbr }^{k_{n-1}+2} \rangle\nonumber \\
 &&+\langle  e_{\lbr b\rbr }e_{\lbr a^i \rbr }e_{\lbr b\rbr }\rangle z(\lbr b \rbr) \langle \cdots e_{\lbr a^{i}\rbr }^{k_{i}+1}   \cdots e_{\lbr b\rbr }^{k_{n-1}+1} e_{\lbr ab\rbr }^{k_{n}+1}\rangle +\text{Length}(<m)\nonumber \\
&=&\text{Length}(<m)+\left\{ \begin{array}{ll}
4\lambda\langle \cdots e_{\lbr a^{i}\rbr }^{k_{i}+1}   \cdots e_{\lbr b\rbr }^{k_{n-1}+1} e_{\lbr ab\rbr }^{k_{n}+1}\rangle, & \text{if}\quad 2|i\vspace{0.2cm}\\
2\lambda\langle \cdots e_{\lbr a^{i}\rbr }^{k_{i}+1}   \cdots e_{\lbr b\rbr }^{k_{n-1}+2} e_{\lbr ab\rbr }^{k_{n}}\rangle+2\langle \cdots e_{\lbr a^{i}\rbr }^{k_{i}+1}   \cdots e_{\lbr b\rbr }^{k_{n-1}} e_{\lbr ab\rbr }^{k_{n}+2}\rangle, & \text{if}\quad 2\nmid i
\end{array}\right.\vspace{0.2cm}\nonumber\\
&=&\text{Length}(<m)+\left\{ \begin{array}{ll}
4\lambda\langle \cdots e_{\lbr a^{i}\rbr }^{k_{i}+1}   \cdots e_{\lbr b\rbr }^{k+1} e_{\lbr ab\rbr }^{k+1}\rangle, & \text{if}\quad 2|i \vspace{0.2cm}\\
4\lambda\langle \cdots e_{\lbr a^{i}\rbr }^{k_{i}+1}   \cdots e_{\lbr b\rbr }^{k+2} e_{\lbr ab\rbr }^{k}\rangle, & \text{if}\quad 2\nmid i
\end{array}\right. .
\eea
In the last equality we have taken $k_{n-1}=k_{n}=k$.
Now we fix $k_1,\cdots, k_{n-2},k$ and set for $1\leq i\leq n-2 $
\bea\label{4}
u_i=\left\{\begin{array}{ll}
\langle \cdots e_{\lbr a^{i}\rbr }^{k_{i}+1}   \cdots e_{\lbr b\rbr }^{k+1} e_{\lbr ab\rbr }^{k+1}\rangle, & \text{if}\quad 2|i , \vspace{0.2cm}\\
\langle \cdots e_{\lbr a^{i}\rbr }^{k_{i}+1}   \cdots e_{\lbr b\rbr }^{k+2} e_{\lbr ab\rbr }^{k}\rangle, & \text{if}\quad 2\nmid i.
\end{array}\right.
\eea
Thus by the equality (\ref{7})=(\ref{8}) from the WDVV equation (\ref {2}) we obtain $n-3$ linear equations of $u_1,\cdots, u_{n-2}$. Consider also the WDVV
\bea\label{3}
&&\langle e_{\lbr a\rbr }^{k_1}\cdots  e_{\lbr a^{n-2}\rbr }^{k_{n-2}} e_{\lbr b\rbr }^{k_{n-1}} e_{\lbr ab\rbr }^{k_{n}}( e_{\lbr a^{n-2}\rbr } e_{\lbr a^{n-2}\rbr }| e_{\lbr b\rbr } e_{\lbr ab\rbr })\rangle\nonumber\\
&=&\langle e_{\lbr a\rbr }^{k_1}\cdots  e_{\lbr a^{n-2}\rbr }^{k_{n-2}} e_{\lbr b\rbr }^{k_{n-1}} e_{\lbr ab\rbr }^{k_{n}}( e_{\lbr a^{n-2}\rbr } e_{\lbr b\rbr }| e_{\lbr a^{n-2}\rbr } e_{\lbr ab\rbr })\rangle.
\eea
Expanding both sides, noting that there exists no nonzero length correlator of the form $\langle  e_{\lbr a^{l}\rbr }e_{\lbr a^{n-2}\rbr }e_{\lbr a^{n-2}\rbr }\rangle$, where $1\leq l\leq n-2$, we obtain
\bea\label{66}
&&\sum_{1\leq l\leq n-2, 2\nmid l}\langle \cdots e_{\lbr a^{l}\rbr }^{k_{l}+1}\cdots e_{\lbr a^{n-2}\rbr }^{k_{n-2}+2}\cdots \rangle \cdot (n-2)\lambda
+\text{Length}(<m)\nonumber\\
&=& 2\lambda\langle \cdots e_{\lbr a^{n-2}\rbr }^{k_{n-2}+1}   \cdots e_{\lbr b\rbr }^{k_{n-1}+1} e_{\lbr ab\rbr }^{k_{n}+1}\rangle+\text{Length}(<m).
\eea
Together with the proceeding $n-3$ equations, we obtain a linear system of $u_1,\cdots,u_{n-2}$. We need to show that this system is nonsingular.
For this, consider the fate of each element $i$ in the set $\{1,\cdots,n-3\}$ under the iterations of the map $\phi$. There are four types:\\
Type 1. After a finite iteration of $\phi$, $i$ maps to $n-2$.\\
Type 2. $i$ is fixed by $\phi$.\\
Type 3. There exists some $r>1$ such that $i=\phi^{(r)}(i)$, i.e., $i,\phi (i),..., \phi^{(r-1)}(i)$ form a cycle under the iterations of $\phi$. We call such $i$ a $\phi$-cyclic element.\\
Type 4. After a finite iteration of $\phi$, $i$ maps to a $\phi$-cyclic element, or a $\phi$-fixed element.\\

For $i$ of type 1, $u_i$ obviously can be determined by the correlators of length less than $m$ and the correlators with less insertions of $e_{\lbr b\rbr}$ and $e_{\lbr ab\rbr}$, since $u_{n-2}$ does, by (\ref{66}). Type 2 is also trivial. \\

For $i$ of type 3, note that the  equations involving $u_{i},\cdots,u_{\phi^{(r-1)}(i)}$ form a subsystem, for which the corresponding matrix of coefficients can be arranged to the form
$$
\left (\begin{array}{ccccc}
4& -1 & & &\\
& 4 & -1 & &  \\
&& \cdots & \cdots & \\
&&& 4 & -1 \\
-1 &&&& 4 \\
\end{array}
\right),
$$
which is easily seen to be nonsingular. Finally, the type 4 cases are reduced to the first three cases. Thus we complete the induction step for the correlators in the form of (\ref{4}).\\

 Again we fix $k_1,\cdots, k_{n-2}, k$ and set for $\leq i\leq n-2 $
\bea\label{5}
v_i=\left\{\begin{array}{ll}
\langle \cdots e_{\lbr a^{i}\rbr }^{k_{i}+1}   \cdots e_{\lbr b\rbr }^{k+1} e_{\lbr ab\rbr }^{k+1}\rangle, & \text{if}\quad 2\nmid i , \vspace{0.2cm}\\
\langle \cdots e_{\lbr a^{i}\rbr }^{k_{i}+1}   \cdots e_{\lbr b\rbr }^{k+2} e_{\lbr ab\rbr }^{k}\rangle, & \text{if}\quad 2| i.
\end{array}\right.
\eea
Consider the WDVV equation
\bea\label{6}
&&\langle e_{\lbr a\rbr }^{k_1}\cdots  e_{\lbr a^{n-2}\rbr }^{k_{n-2}} e_{\lbr b\rbr }^{k} e_{\lbr ab\rbr }^{k}( e_{\lbr a^{i}\rbr } e_{\lbr a^{i}\rbr }| e_{\lbr ab\rbr } e_{\lbr ab\rbr })\rangle\nonumber\\
&=&\langle e_{\lbr a\rbr }^{k_1}\cdots  e_{\lbr a^{n-2}\rbr }^{k_{n-2}} e_{\lbr b\rbr }^{k} e_{\lbr ab\rbr }^{k}( e_{\lbr a^{i}\rbr } e_{\lbr ab\rbr }| e_{\lbr a^{i}\rbr } e_{\lbr ab\rbr })\rangle,
\eea
where $1\leq i\leq n-2$. Thus in a similar way we have
\bea\label{9}
&&\sum_{1\leq l\leq n-2, 2|l}\langle \cdots e_{\lbr a^{i}\rbr }^{k_{i}+2}\cdots e_{\lbr a^{l}\rbr }^{k_{l}+1}\cdots \rangle \cdot (n-2)\lambda
+\langle \cdots  e_{\lbr a^{\phi(i)}\rbr }^{k_{\phi(i)}+1} \cdots e_{\lbr b\rbr }^{k} e_{\lbr ab\rbr }^{k+2}\rangle\nonumber\\
&&\cdot\big(1\cdot(1-\delta_{i,\frac{n-2}{2}})+2\delta_{i,\frac{n-2}{2}}\big)\lambda+\text{Length}(<m)=4v_i,
\eea
and
\bea\label{10}
&&\sum_{1\leq l\leq n-2, 2| l}\langle \cdots e_{\lbr a^{l}\rbr }^{k_{l}+1}\cdots e_{\lbr a^{n-2}\rbr }^{k_{n-2}+2}\cdots \rangle \cdot (n-2)\lambda
+\text{Length}(<m)\nonumber\\
&=& 2\lambda\langle \cdots e_{\lbr a^{n-2}\rbr }^{k_{n-2}+1}   \cdots e_{\lbr b\rbr }^{k} e_{\lbr ab\rbr }^{k+2}\rangle+\text{Length}(<m).
\eea
The argument for determination of $u_i$'s applies for $v_i$'s. \\

 Thus the proof of the lemma when $2|n$ is completed.\\

When $2\nmid n$, we still set $u_i$ as in (\ref{4}). We consider the WDVV equation (\ref{2}) and (\ref{3}), and one easily checks that (\ref{7}) (\ref{8}) still holds, while (\ref{66}) turns into
\bea\label{67}
&&\sum_{1\leq l\leq n-2, 2\nmid l}\langle \cdots e_{\lbr a^{l}\rbr }^{k_{l}+1}\cdots e_{\lbr a^{n-2}\rbr }^{k_{n-2}+2}\cdots \rangle \cdot (n-2)\lambda
+\text{Length}(<m)\nonumber\\
&=& \lambda\langle \cdots e_{\lbr a^{n-2}\rbr }^{k_{n-2}+1}   \cdots e_{\lbr b\rbr }^{k+2} e_{\lbr ab\rbr }^{k}\rangle+\lambda\langle \cdots e_{\lbr a^{n-2}\rbr }^{k_{n-2}+1}   \cdots e_{\lbr b\rbr }^{k} e_{\lbr ab\rbr }^{k+2}\rangle+\text{Length}(<m).\nonumber\\
\eea
The argument for the nonsingularity of the linear system still holds. A similar argument holds for the correlator of the form (\ref{5}).
\hfill\qedsymbol

\begin{lemma}\label{10}
The length $m(\geq 4)$ correlators without insertions of the form $e_{\lbr a^k \rbr}$, $1\leq k\leq n-2$, are uniquely determined by the correlators of length less than $m$ and the correlators with at least one insertion of the form $e_{\lbr a^i \rbr}$, $1\leq i\leq n-2$, and the correlators of the form  $\langle e_{\lbr b \rbr}^k \rangle$ or $\langle e_{\lbr ab \rbr}^k \rangle$ and the WDVV equations.
\end{lemma}
Proof:
Note that $\langle e_{\lbr b \rbr}^{k_{n-1}}  e_{\lbr ab \rbr}^{k_{n}} \rangle\neq 0$ implies that $k_{n-1}$ and $k_n$ are both even. Therefore by (\ref{11}) and (\ref{12}) we can inductively determines $\langle e_{\lbr b \rbr}^{k_{n-1}}  e_{\lbr ab \rbr}^{k_{n}} \rangle$ from correlators $\langle e_{\lbr b \rbr}^{k_{n-1}+k_{n} }\rangle$ or $\langle e_{\lbr ab \rbr}^{k_{n-1}+k_{n}} \rangle$, and correlators with at least one insertion $e_{\lbr a^{i} \rbr}$, $1\leq i\leq n-2$, and the correlators of length less than $m$ . \hfill\qedsymbol
\vspace{0.5cm}\\
Combining lemma \ref{1} and lemma \ref{10}, by induction on the length of the correlators we obtain
\begin{theorem}\label{13}
The correlators $\langle \cdot \rangle ^{[\mathbb{C}^2/\hat{D}_{n}]}$ of length at least four are determined by correlators of length three, the correlators with only insertions of the form $e_{\lbr a^i \rbr}$, $1\leq i\leq n-2$, the correlators $\langle e_{\lbr b \rbr}^k \rangle$, $\langle e_{\lbr ab \rbr}^k \rangle$, $k\geq 4$, and the WDVV equations.\hfill\qedsymbol
\end{theorem}

\begin{corollary}\label{30}
Conjecture \ref{27} holds for $\hat{D}_{n}$ if $n=2^{m}+2$ for some positive integer $m\geq 1$.
\end{corollary}
Proof: It is easily seen that Conjecture \ref{27} holds for length three correlators.  The right handside of (\ref{25}) satisfies the WDVV equation automatically, since it is the genus 0 primary Gromov-Witten invariants of the resolution of $\mathbb{C}^2/\hat{D}_{n}$. By theorem \ref{13}, it suffices to verify (\ref{25}) for the special correlators in theorem \ref{13}. By the results of compatibility of section 4, and theorem \ref{31},  it reduces to verify this for $\hat{D}_{4}$. But $\hat{D}_{4}$ has an automorphism sending $a$ to $b$, and the pullback of the representation $\rho_{1}$ is equivalent to itself. Therefor we complete the proof.\hfill\qedsymbol

\begin{theorem}\label{32}
Conjecture \ref{27} holds for $\hat{D}_{n}$ when $n\geq 4$.
\end{theorem}
Proof: When $n$ is even, the statement follows from Corollary \ref{30} and the polynomiality theorem \ref{70}. When $n$ is odd, by propsition \ref{105} we have
\ben
\langle e_{\lbr b\rbr}^{2k+2}\rangle^{[\mathbb{C}^2/\hat{D}_{2n-2}]}=\frac{1}{2}\sum_{k_b=0}^{2k+2} \binom{2k+2}{k_b}\langle e_{\lbr b\rbr}^{k_b}e_{\lbr ab\rbr}^{2k+2-k_b}\rangle^{[\mathbb{C}^2/\hat{D}_{n}]}.
\een
The statement for $\hat{D}_{2n-2}$ will impose another constraint to the correlators of $\hat{D}_{n}$ by the compatibility proved in section 4.2. Note that for a given $k$, exactly one of $\langle e_{\lbr b\rbr}^{k+1}e_{\lbr ab\rbr}^{k+1}\rangle$ and $\langle e_{\lbr b\rbr}^{k}e_{\lbr ab\rbr}^{k+2}\rangle$ is nonzero, according to the parity of $k$. This consideration together with the equation (\ref{11}) and (\ref{12}) completely determines correlators $\langle e_{\lbr b \rbr}^{2k+2} \rangle$ and $\langle e_{\lbr ab \rbr}^{2k+2} \rangle$ of length $2k+2$ from the correlators of length less than $2k+2$ and the correlators with at least one insertion of the form $e_{\lbr a^i \rbr}$, $1\leq i\leq n-2$ . Therefore applications of lemma \ref{1}, lemma \ref{10}, the WDVV equations give the result.
\hfill\qedsymbol

\begin{theorem}\label{107}
Conjecture \ref{28} holds for ${D}_{n}$ when $n\geq 4$.
\end{theorem}
Sketch of the proof: The WDVV induction is treated in subsection 5.1. The induction from normal subgroups and the polynomiality is very similar to those for  $\hat{D}_{n}$, and we omit them. The induction starts from $D_4$, which is an abelian group and the corresponding conjecture has  been proved in \cite{BGh0708}. \hfill\qedsymbol

\section{Some combinatorics of fractional Bernoulli numbers}
\subsection{Properties of $\frac{1}{\zeta_{n}^{b}-1}$}
For a natural number $n$, denote $\zeta_{n}=\exp (\frac{2\pi i}{n})$, and define $$\alpha_{n,b}=\frac{1}{\zeta_{n}^{b}-1},$$ for $b=1,\cdots, n$. We formally set
\bea
\alpha_{n,0}=0.
\eea
Let $p_m$ be the $m$-th Newton symmetric function, i.e., $p_m(x_1,\cdots)=\sum_{i}x_{i}^{m}$ for variables $x_1,\cdots$. Then we have
\begin{lemma}\label{15}
For every fixed $m\geq 0$, $p_m (\alpha_{n,1},\cdots,\alpha_{n,n-1})$ is a polynomial of $n$, when $n$ varies in $2,3,\cdots$. Consequently, any fixed symmetric function valued at  $\alpha_{n,1},\cdots,\alpha_{n,n-1}$ depends  polynomially on $n$.
\end{lemma}
Proof: We compute the generating function
\ben
\sum_{k=0}^{\infty}\sum_{b=1}^{n-1}\frac{z^k}{(\zeta_{n}^{b}-1)^k}&=& \sum_{b=1}^{n-1}\frac{1}{1-{\frac{z}{\zeta_{n}^{b}-1}}}\\
&=& n-1+ \sum_{b=1}^{n-1}\frac{z}{\zeta_{n}^{b}-1-z}\\
&=&n+\sum_{b=0}^{n-1}\frac{z}{\zeta_{n}^{b}-1-z}\\
&=& n-z\frac{d}{dz}\log \prod_{b=0}^{n-1}(1+z-\zeta_{n}^{b})\\
&=&  n-z\frac{d}{dz}\log\big((1+z)^{n}-1\big)\\
&=& n+\frac{n z(1+z)^{n-1}}{1-(z+1)^n}\\
&=& (n-1)+\frac{n-1}{2}z+\frac{-5 + 6 n - n^2}{12}z^2+\frac{3 - 4 n + n^2}{8}z^3\\
 && +\frac{-251 + 360 n - 110 n^2 + n^4}{720}z^4+\frac{95 - 144 n + 50 n^2 - n^4}{288}z^5+\cdots
\een
In general, it is easily seen that the coefficient of $z^k$ in the expansion of $\frac{z}{1-(1+z)^n}$ at $z=0$ is a polynomial in $n$, for any fixed $k$.\hfill\qedsymbol

The Bernoulli polynomials $B_m(x)$ are defined by
\bea
\frac{te^{tx}}{e^{t}-1}=\sum_{m=0}^{\infty}\frac{B_m(x)t^m}{m!},
\eea
in particular, $B_m=B_m(0)$ are the ordinary Bernoulli numbers.
Define
\bea
\beta_{m_{1},\cdots,m_{k}}(n,b)=\sum_{a=0}^{n-1}\zeta_{n}^{ab}\prod_{j=1}^{k}B_{m_{j}} (\frac{a}{n}),
\eea
for $b\in\mathbb{Z}$, thus
\bea
\sum_{b=0}^{n-1}\beta_{m_{1},\cdots,m_{k}}(n,b)=n\prod_{j=1}^{k}B_{m_{j}}.
\eea
Sometimes we use capital letters to stand for  sets of subscripts for convenience, i.e., $\beta_{M}(n,b):=\beta_{m_{1},\cdots,m_{k}}(n,b)$ when $M=\{m_{1},\cdots,m_{k}\}$. One easily checks that
\bea
\sum_{i=0}^{n-1}\beta_{M}(n,a+i)\beta_{N}(n,b-i)=n\beta_{M\coprod N}(n,a+b),
\eea
and
\bea\label{50}
\beta_{M}(2n,b)+\beta_{M}(2n,n+b)=2\beta_{M}(n,b).
\eea

Consider the generating function, we have
\bea\label{14}
&&\sum_{m_{1}=0}^{\infty}\cdots\sum_{m_{k}=0}^{\infty}\frac{\beta_{m_{1},\cdots,m_{k}}(n,b)t_{1}^{m_{1}}\cdots t_{k}^{m_{k}}}{m_{1}!\cdots m_{k}!}\nonumber\\
&=&\sum_{a=0}^{n-1}\zeta_{n}^{ab}\prod_{j=1}^{k}\frac{t_{j}e^{\frac{at}{n}}}{e^{t_j}-1}\nonumber\\
&=&\prod_{j=1}^{k}\frac{t_{j}}{e^{t_j}-1}\cdot \frac{\prod_{j=1}^{k}e^{t_j}-1}{\zeta_{n}^{b}\prod_{j=1}^{k}e^{\frac{t_j}{n}}-1}.
\eea

\begin{lemma}\label{23}
There exists a function $f_{m_1,\cdots,m_k}(x,y)=y^{m-1}f_{m-1}(x)+y^{m-2}f_{m-2}(x)+\cdots+f_{0}(x)$, where $m=\sum_{i=1}^{k}m_i$, and $f_{i}(x)$ depending only on $m_1,\cdots,m_k$,  is a polynomial of $x$ for $i=0,\cdots, m-1$,  such that $$\beta_{m_1,\cdots,m_k}(n,b)=f_{m_1,\cdots,m_k} (\alpha_{n,b},\frac{1}{n}),$$
for all $n\geq 1$, $1\leq b\leq n-1$.
\end{lemma}
Proof:  For $b=1,\cdots, n-1$, this follows from (\ref{14}) and the expansion
\ben
\frac{1}{\zeta_{n}^{b}e^{\frac{t}{n}}-1}=\frac{1}{\zeta_{n}^{b}-1} \sum_{k=0}^{\infty}(-1)^{k}(\frac{1}{\zeta_{n}^{b}-1}+1)^{k}(e^{\frac{t}{n}}-1)^{k}.
\een
\hfill\qedsymbol

We treat $\alpha_{n,0}$ separately. Take $k=1$ in (\ref{14}) we obtain
\bea\label{45}
\beta_{m}(n,0)&=&\frac{B_m}{n^{m-1}},
\eea

\bea\label{42}
\beta_{m}(nc,bc)&=&\frac{\beta_{m}(n,b)}{c^{m-1}},
\eea
where $c\in \mathbb{Z}$.
In general we have the
\begin{lemma}\label{22}
For fixed $m_{1},\cdots,m_{k}$ and fixed $n$, $b$, $\beta_{m_{1},\cdots,m_{k}}(nc,bc)$ as a function of the positive integer $c$ is a Laurent polynomial of $c$. In particular, $\beta_{m_{1},\cdots,m_{k}}(c,0)$ and $\beta_{m_{1},\cdots,m_{k}}(2c,c)$ is a Laurent polynomials  of $c$.
\end{lemma}
Proof: This is obvious by expanding
\ben
&&\sum_{m_{1}=0}^{\infty}\cdots\sum_{m_{k}=0}^{\infty}\frac{\beta_{m_{1},\cdots,m_{k}}(nc,bc)(ct_{1}^{m_{1}})\cdots (ct_{k}^{m_{k}})}{m_{1}!\cdots m_{k}!}\\
&=& \frac{c^k t_1\cdots t_k}{(e^{ct_1}-1)\cdots(e^{ct_k}-1)}\cdot \frac{e^{c(t_1+\cdots+t_k)}-1}{\zeta_{nc}^{bc}e^{\frac{t_1+\cdots t_k}{n}}-1}\\
&=&  \frac{c^k (e^{(c-1)\sum_{i=1}^{k}t_i}+\cdots+e^{\sum_{i=1}^{k}t_i}+1)}{\prod_{i=1}^{k}(e^{(c-1)t_i}+\cdots+e^{t_i}+1)}\cdot \prod_{j=1}^{k}\frac{t_{j}}{e^{t_j}-1}\cdot \frac{\prod_{j=1}^{k}e^{t_j}-1}{\zeta_{n}^{b}\prod_{j=1}^{k}e^{\frac{t_j}{n}}-1}\\
&=& \frac{c^k (e^{(c-1)\sum_{i=1}^{k}t_i}+\cdots+e^{\sum_{i=1}^{k}t_i}+1)}{\prod_{i=1}^{k}(e^{(c-1)t_i}+\cdots+e^{t_i}+1)}\cdot
\sum_{m_{1}=0}^{\infty}\cdots\sum_{m_{k}=0}^{\infty}\frac{\beta_{m_{1},\cdots,m_{k}}(n,b)t_{1}^{m_{1}}\cdots t_{k}^{m_{k}}}{m_{1}!\cdots m_{k}!}.
\een
\hfill\qedsymbol\vspace{1cm}\\

\subsection{Summing fractional Bernoulli numbers over a Feymann diagram}
\begin{definition}
Let $\Gamma=\{V(\Gamma),E(\Gamma)\}$ be a connected graph, with each edge decorated by a finite set consisting of positive integers. We call $\Gamma$ with such a decoration a \textit{decorated Feymann diagram}.
\end{definition}
For a given integer $n\geq 1$, and a decorated Feymann diagram, we associate a number $S_{n}(\Gamma)$ as follows. For the purpose of this article, we focus on tree diagrams. For a vertex $v\in V(\Gamma)$, let $E(v)$ be the edges associated to $v$. For an edge $e\in E(\Gamma)$, denote the two endpoints of $e$ by $u(e)$ and $w(e)$, with arbitrary choice between the two. We say that $u(e)$ and $w(e)$ are neighboring. \\

We shall also need a filtration of $\Gamma$. Let $\Gamma^{(0)}=\Gamma$. Cut out the vertices of valence 1 (i.e. \emph{leaves}) in $\Gamma$ and the edges incident to them, we obtain a subtree, which we denote by $\Gamma^{(1)}$. We continue this procedure to obtain a finite filtration $\Gamma^{(0)}\supset \Gamma^{(1)}\supset \cdots \supset \Gamma^{(N+1)}=\emptyset$. Note that when $\Gamma^{(k+1)}\neq \emptyset$, any two vertices in $V(\Gamma^{(k)})\backslash V(\Gamma^{(k+1)})$ is not connected by any edge in $\Gamma$.\\

First, for an edge $e$ decorated by $(m_{1,e},\cdots,m_{k,e})$, we associate the set of  products $P_{e}=\{B_{m_{1,e}}(\frac{{a}_{e}}{n})\cdots B_{m_{k,e}}(\frac{{a}_{e}}{n})|0\leq a_{e}\leq n-1\}$. Next we sum all the products of the form
\ben
\prod_{e\in E(\Gamma)}B_{m_{1,e}}(\frac{{a}_{e}}{n})\cdots B_{m_{k,e}}(\frac{{a}_{e}}{n}),
\een
where each $a_{e}$ running over $0$ to $n-1$, under the restriction that, for each vertex $v\in V(\Gamma)$ of valence $\geq 2$, we demand that $\sum_{e\in E(v)}a_{e}\equiv0(\text{mod}\hspace{0.1cm} n)$. In other words, every vertex $v$ contributes a $\delta$-function $\delta(\sum_{e\in E(v)}a_{e})$ in the $\text{mod}\hspace{0.1cm} n$ sense.
\begin{example}
Let $\Gamma_{1}(m)$ be the diagram
$$ \xy
(0,0); (10,0), **@{-};
(0,0)*+{\bullet};(10,0)*+{\bullet};(5,2)*+{m};
\endxy ,
$$
then
\ben
S_{n}(\Gamma_{1})=\sum_{a=0}^{n-1}B_{m}(\frac{a}{n})=\beta_{m}(n,0),
\een
which is a Laurent polynomial of $n$ by (\ref{45}).
\end{example}

\begin{example}\label{46}
Let $\Gamma_{2}(m_1, m_2, m_3)$ be the diagram
$$ 
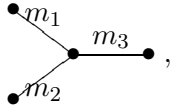
\xy
(0,0); (10,0), **@{-};
(0,0); (-8,6), **@{-};(0,0); (-8,-6), **@{-};(-8,-6)*+{\bullet};(-4,-5)*+{m_2};(-4,5)*+{m_1};(-8,6)*+{\bullet};
(10,0)*+{\bullet};(0,0)*+{\bullet};(5,2)*+{m_3};
\endxy,
$$
then
\ben
S_{n}(\Gamma_{2})&=&\sum_{0\leq a_{i}\leq n-1, i=1,2,3, a_1+a_2+a_3\equiv 0(\text{mod}\hspace{0.1cm} n)}B_{m_1}(\frac{a_1}{n})B_{m_2}(\frac{a_2}{n})B_{m_3}(\frac{a_3}{n})\\
&=& \sum_{a_{1}=0}^{n-1}\sum_{a_{2}=0}^{n-1}\sum_{a_{3}=0}^{n-1}\sum_{b=0}^{n-1}
\frac{\zeta_{n}^{b(a_1+a_2+a_3)}}{n}B_{m_1}(\frac{a_1}{n})B_{m_2}(\frac{a_2}{n})B_{m_3}(\frac{a_3}{n})\\
&=&\frac{1}{n}\sum_{b=0}^{n-1}\beta_{m_1}(n,b)\beta_{m_2}(n,b)\beta_{m_3}(n,b),
\een
which is a Laurent polynomial of $n$, by lemma \ref{15} and lemma \ref{23}. The same proof applies  for all \emph{star} diagrams as well.
\end{example}

\begin{example}\label{106}
Let $\Gamma_{3}(m_1,\cdots,m_7)$ be the diagram
$$ 
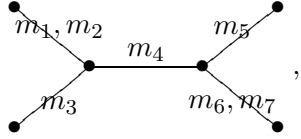
\xy
(0,0); (15,0), **@{-};
(0,0); (-10,8), **@{-};(0,0); (-10,-8), **@{-};(-10,-8)*+{\bullet};(-4,-5)*+{m_3};(-4,5)*+{m_1,m_2};(-10,8)*+{\bullet};
(15,0)*+{\bullet};(0,0)*+{\bullet};(7.5,2)*+{m_4};
(15,0); (25,8), **@{-};(15,0); (25,-8), **@{-};(25,-8)*+{\bullet};(19,-5)*+{m_6,m_7};(19,5)*+{m_5};(25,8)*+{\bullet};
\endxy,
$$
then
\ben
S_{n}(\Gamma_{3})&=&\sum_{0\leq a_{i}\leq n-1, i=1,\cdots,5, a_1+a_2+a_3\equiv a_3+a_4+a_5\equiv 0(\text{mod}\hspace{0.1cm}n)}
B_{m_1}(\frac{a_1}{n})B_{m_2}(\frac{a_1}{n})B_{m_3}(\frac{a_2}{n})\\
&&B_{m_4}(\frac{a_3}{n})B_{m_5}(\frac{a_4}{n})
\cdot B_{m_6}(\frac{a_5}{n})B_{m_7}(\frac{a_5}{n})\\
&=&  \sum_{a_{1}=0}^{n-1}\cdots\sum_{a_{5}=0}^{n-1}\sum_{b_1=0}^{n-1}\sum_{b_2=0}^{n-1}
\frac{\zeta_{n}^{b_{1}(a_1+a_2+a_3)}}{n}\frac{\zeta_{n}^{b_{2}(a_3+a_4+a_5)}}{n}\\
&&B_{m_1}(\frac{a_1}{n})B_{m_2}(\frac{a_1}{n})B_{m_3}(\frac{a_2}{n})B_{m_4}(\frac{a_3}{n})B_{m_5}(\frac{a_4}{n})\\
&=&\frac{1}{n^2}\sum_{b_1=0}^{n-1}\sum_{b_2=0}^{n-1}\beta_{m_1,m_2}(n,b_1)\beta_{m_3}(n,b_1)\beta_{m_4}(n,b_1+b_2)
\beta_{m_5}(n,b_2)\beta_{m_6,m_7}(n,b_2).
\een
We want to show that $S_{n}(\Gamma_{3})$ is a Laurent polynomial of $n$. By the last equality and  lemma \ref{23}, it suffices to show
\ben
\sum_{b_1=0}^{n-1}\sum_{b_2=0}^{n-1}(\alpha_{n,b_1})^{k_1}(\alpha_{n,b_2})^{k_2}(\alpha_{n,b_1+b_2})^{k_3}
\een
is a Laurent polynomial of $n$, for any fixed $k_1,k_2,k_3\geq 0$. Consider the generating function
\ben
&&\sum_{k_1=0}^{\infty}\sum_{k_2=0}^{\infty}\sum_{k_3=0}^{\infty}
\sum_{b_1=0}^{n-1}\sum_{b_2=0}^{n-1}(\alpha_{n,b_1})^{k_1}(\alpha_{n,b_2})^{k_2}(\alpha_{n,b_1+b_2})^{k_3}t_{1}^{k_1}t_{2}^{k_2}t_{3}^{k_3}\\
&=&
\sum_{b_1=0}^{n-1}\sum_{b_2=0}^{n-1}\frac{\zeta_{n}^{b_1}-1}{\zeta_{n}^{b_1}-1-t_1}\frac{\zeta_{n}^{b_2}-1}{\zeta_{n}^{b_2}-1-t_2}
\frac{\zeta_{n}^{b_1+b_2}-1}{\zeta_{n}^{b_1+b_2}-1-t_3}\\
&=& \sum_{b_1=0}^{n-1}\sum_{b_2=0}^{n-1}\sum_{b_3=0}^{n-1}\sum_{c=0}^{n-1}\frac{\zeta_{n}^{c(b_1+b_2-b_3)}}{n}
\frac{\zeta_{n}^{b_1}-1}{\zeta_{n}^{b_1}-1-t_1}\frac{\zeta_{n}^{b_2}-1}{\zeta_{n}^{b_2}-1-t_2}
\frac{\zeta_{n}^{b_3}-1}{\zeta_{n}^{b_3}-1-t_3}\\
&=&\frac{1}{n}\sum_{b_2=0}^{n-1}\sum_{c=0}^{n-1}\frac{\zeta_{n}^{b_2(c+1)}-\zeta_{n}^{b_2 c}}{\zeta_{n}^{b_2}-1-t_2}\sum_{b_1=0}^{n-1}\frac{\zeta_{n}^{b_1(c+1)}-\zeta_{n}^{b_1 c}}{\zeta_{n}^{b_1}-1-t_1}
\sum_{b_3=0}^{n-1}\frac{\zeta_{n}^{b_3(1-c)}-\zeta_{n}^{-b_3 c}}{\zeta_{n}^{b_3}-1-t_3}\\
&=& \frac{1}{n}\sum_{b_2=0}^{n-1}\sum_{c=1}^{n}\frac{\zeta_{n}^{b_2(c+1)}-\zeta_{n}^{b_2 c}}{\zeta_{n}^{b_2}-1-t_2}\Big(
n\delta_{c,n}+\frac{nt_1 (1+t_1)^{c-1}}{1-(1+t_1)^{n}}\Big)\Big(
n\delta_{c,n}+\frac{nt_3 (1+t_3)^{n-c}}{1-(1+t_3)^{n}}\Big).
\een
Here the last equality follows inductively from
\ben
\sum_{b=0}^{n-1}\frac{\zeta_{n}^{b (c+1)}}{\zeta_{n}^{b}-1-t}=n\delta_{0,c}+(1+t)\sum_{b=0}^{n-1}\frac{\zeta_{n}^{b c}}{\zeta_{n}^{b}-1-t}
\een
for $0\leq c\leq n-1$, and
\ben
\sum_{b=0}^{n-1}\frac{1}{\zeta_{n}^{b}-1-t}=\frac{n(1+t)^{n-1}}{1-(1+t)^{n}}
\een
as in the proof of lemma \ref{15}.\\
For $1\leq b_2\leq n-1$,
\ben
\sum_{c=1}^{n}(\frac{1+t_1}{1+t_3})^{c}\zeta_{n}^{cb_2}=-\big(1-(\frac{1+t_1}{1+t_3})^{n}\big)(1+\alpha_{n,b_2})\cdot
\frac{1+t_1}{1+t_1 -\alpha_{n,b_2}(t_3-t_1)},
\een
thus the statement reduces to the case of example \ref{46}. This proof applies for any $\Gamma$ with $\Gamma^{(2)}=0$ as well.
\end{example}

\begin{proposition}\label{43}
For any decorated \emph{tree} Feymann diagram $\Gamma$,
\ben
\sum_{0\leq b_v\leq n-1, v\in V(\Gamma^{1})}\Bigg(\sum_{v\in V(\Gamma^{1})}(\alpha_{n, b_v})^{k_v}\sum_{e\in E(\Gamma^{1})}(\alpha_{n, b_{u(e)}+b_{w(e)}})^{k_e}\Bigg)
\een
is a Laurent polynomial of $n$, with fixed natural numbers $k_v$'s and $k_e$'s.
\end{proposition}
Proof: We prove this proposition by induction along the stratification of $\Gamma$. For $v_1\in V(\Gamma^{(k)})\backslash V(\Gamma^{(k+1)})$, and $v_2\in V(\Gamma^{(k+1)})$ such that there is an edge $e\in E(\Gamma)$ connecting $v_1$ and $v_2$, consider the generating function
\ben
&&\sum_{k_1=0}^{\infty}\sum_{k_2=0}^{\infty}\sum_{k_e=0}^{\infty}
\sum_{b_1=0}^{n-1}(\alpha_{n,b_1})^{k_1}(\alpha_{n,b_2})^{k_2}(\alpha_{n,b_1+b_2})^{k_e}t_{1}^{k_1}t_{2}^{k_2}t_{e}^{k_e}\\
&=&
\sum_{b_1=0}^{n-1}\frac{\zeta_{n}^{b_1}-1}{\zeta_{n}^{b_1}-1-t_1}\frac{\zeta_{n}^{b_2}-1}{\zeta_{n}^{b_2}-1-t_2}
\frac{\zeta_{n}^{b_1+b_2}-1}{\zeta_{n}^{b_1+b_2}-1-t_e}\\
&=& \sum_{b_1=0}^{n-1}\sum_{b_e=0}^{n-1}\sum_{c=0}^{n-1}\frac{\zeta_{n}^{c(b_1+b_2-b_e)}}{n}
\frac{\zeta_{n}^{b_1}-1}{\zeta_{n}^{b_1}-1-t_1}\frac{\zeta_{n}^{b_2}-1}{\zeta_{n}^{b_2}-1-t_2}
\frac{\zeta_{n}^{b_e}-1}{\zeta_{n}^{b_e}-1-t_e},
\een
thus the same argument as in example \ref{106} reduce $\Gamma^{(k)}$ to the diagram $\Gamma^{(k)}\backslash\{v_1,e\}$. This process finally terminates at a star diagram since the diagram is tree.
\hfill\qedsymbol
\\

\begin{theorem}\label{36}
For a decorated Feymann \emph{tree} diagram $\Gamma$, $S_{n}(\Gamma)$ is a Laurent polynomial of $n$.
\end{theorem}
Proof: Since for any $v\in V(\Gamma^{1})$, we have
\ben
&&\sum_{0\leq a_e\leq n-1, e\in E(v), \sum_{e\in E(v)}a_e\equiv 0(\text{mod}\hspace{0.1cm}n)}(\cdot)\\
&=& \frac{1}{n}\sum_{0\leq a_e\leq n-1, e\in E(v)}\sum_{b_v=0}^{n-1}\zeta_{n}^{b_v\sum_{e\in E(v)}a_e}(\cdot),
\een
by lemma \ref{23} and lemma \ref{22}, we are left to proposition \ref{43}.
\hfill\qedsymbol

\section{Polynomiality of $\langle e_{\lbr b \rbr}^{2m} \rangle^{[\mathbb{C}^2 /\hat{D}_n]}$ for even $n$ }
In this section, firstly we write the differential operator $\big(\frac{A_{p+1}(V_{\rho_{1}})z^p}{(p+1)!}\big)^{\wedge}$ from the flat coordinates to the coordinates associated to the semisimple basis of the untwisted Gromov-Witten potential, and then to the mixed coordinates (\ref{108}). Secondly we make some manipulations and then modify it into a more simple operator (\ref{109}), which gives the same correlators we need, by the symmetry of the total descendant potential that the differential operators act on. Finally we make some \emph{d\'{e}vissages} of the differential operators and the trees to prove the (Laurent) polynomiality of the correlators we concern.\\

According to \cite{JK}, the Gromov-Witten potential of $\mathcal{B}G$ takes a simple form in the semisimple basis given by the characters of the irreducible representations of $G$.  We follow the notations of \cite{JK}, and the notations for the representations of $\hat{D}_{n}$ in section 3.
\ben
f_{\psi_{1}}&=& \frac{1}{4n-8}(e_{\lbr 1\rbr}+\sum_{1}^{n-3}e_{\lbr a^k \rbr}+e_{\lbr a^{n-2}\rbr}+e_{\lbr b\rbr}+e_{\lbr ab\rbr}),\\
f_{\psi_{2}}&=& \frac{1}{4n-8}(e_{\lbr 1\rbr}+\sum_{1}^{n-3}e_{\lbr a^k \rbr}+e_{\lbr a^{n-2}\rbr}-e_{\lbr b\rbr}-e_{\lbr ab\rbr}),\\
f_{\psi_{3}}&=& \frac{1}{4n-8}(e_{\lbr 1\rbr}+\sum_{1}^{n-3}(-1)^{k}e_{\lbr a^k \rbr}+e_{\lbr a^{n-2}\rbr}-e_{\lbr b\rbr}+e_{\lbr ab\rbr}),\\
f_{\psi_{4}}&=& \frac{1}{4n-8}(e_{\lbr 1\rbr}+\sum_{1}^{n-3}(-1)^{k}e_{\lbr a^k \rbr}+e_{\lbr a^{n-2}\rbr}+e_{\lbr b\rbr}-e_{\lbr ab\rbr}),\\
f_{\rho_{l}}&=& \frac{2}{4n-8}(2e_{\lbr 1\rbr}+\sum_{1}^{n-3}(\omega^{kl}+\omega^{-kl})e_{\lbr a^k \rbr}+2(-1)^{l}e_{\lbr a^{n-2}\rbr}),\\
\een
for $l=1,\cdots,n-3$. Denote the coordinates according to the basis $\{f_{\alpha}\}$ by $\{u_{j}^{\alpha}\}$, with the subscripts $j$ indicating the degree of powers of the $\psi$-classes, i.e., $\sum_{\alpha}f_{\alpha}u_{j}^{\alpha}=\sum_{\beta}e_{\beta}t_{j}^{\beta}$. Then we have
\ben
t_{j}^{\lbr b\rbr}&=&\frac{1}{4n-8}u_{j}^{\psi_{1}} -\frac{1}{4n-8}u_{j}^{\psi_{2}}-\frac{1}{4n-8}u_{j}^{\psi_{3}} +\frac{1}{4n-8}u_{j}^{\psi_{4}},\\
t_{j}^{\lbr ab\rbr}&=& \frac{1}{4n-8}u_{j}^{\psi_{1}} -\frac{1}{4n-8}u_{j}^{\psi_{2}}+\frac{1}{4n-8}u_{j}^{\psi_{3}} -\frac{1}{4n-8}u_{j}^{\psi_{4}},\\
t_{j}^{\lbr 1\rbr}&=& \frac{1}{4n-8}( u_{j}^{\psi_{1}}+u_{j}^{\psi_{2}}+u_{j}^{\psi_{3}}+u_{j}^{\psi_{4}}+4\sum_{l=1}^{n-3}u_{j}^{\rho_{l}}),\\
t_{j}^{\lbr a^k \rbr}&=& \frac{1}{4n-8}( u_{j}^{\psi_{1}}+u_{j}^{\psi_{2}}+(-1)^{k}u_{j}^{\psi_{3}}+(-1)^{k}u_{j}^{\psi_{4}}+
2\sum_{l=1}^{n-3}(\omega^{kl}+\omega^{-kl})u_{j}^{\rho_{l}}),\\
t_{j}^{\lbr a^{n-2}\rbr}&=& \frac{1}{4n-8}( u_{j}^{\psi_{1}}+u_{j}^{\psi_{2}}+u_{j}^{\psi_{3}}+u_{j}^{\psi_{4}}+4\sum_{l=1}^{n-3}(-1)^{l}u_{j}^{\rho_{l}}),\\
\een
and the inverse transform
\ben
u_{j}^{\psi_1}&=& t_{j}^{\lbr 1\rbr}+ 2\sum_{k=1}^{n-3}t_{j}^{\lbr a^k\rbr}+t_{j}^{\lbr a^{n-2}\rbr}+(n-2)t_{j}^{\lbr b\rbr}+(n-2)t_{j}^{\lbr ab\rbr},\\
u_{j}^{\psi_2}&=&  t_{j}^{\lbr 1\rbr}+ 2\sum_{k=1}^{n-3}t_{j}^{\lbr a^k\rbr}+t_{j}^{\lbr a^{n-2}\rbr}-(n-2)t_{j}^{\lbr b\rbr}-(n-2)t_{j}^{\lbr ab\rbr},\\
u_{j}^{\psi_3}&=&  t_{j}^{\lbr 1\rbr}+ 2\sum_{k=1}^{n-3}(-1)^{k}t_{j}^{\lbr a^k\rbr}+t_{j}^{\lbr a^{n-2}\rbr}-(n-2)t_{j}^{\lbr b\rbr}+(n-2)t_{j}^{\lbr ab\rbr},\\
u_{j}^{\psi_4}&=&  t_{j}^{\lbr 1\rbr}+ 2\sum_{k=1}^{n-3}(-1)^{k}t_{j}^{\lbr a^k\rbr}+t_{j}^{\lbr a^{n-2}\rbr}+(n-2)t_{j}^{\lbr b\rbr}-(n-2)t_{j}^{\lbr ab\rbr},\\
u_{j}^{\rho_l}&=& t_{j}^{\lbr 1\rbr}+\sum_{k=1}^{n-3}(\omega^{lk}+\omega^{-lk})t_{j}^{\lbr a^k\rbr}+(-1)^{l}t_{j}^{\lbr a^{n-2}\rbr}.
\een
In the following we shall also use a \emph{mixed} coordinates,
\bea\label{108}
u_{j}^{\rho_l}&=& t_{j}^{\lbr 1\rbr}+\sum_{k=1}^{n-3}(\omega^{lk}+\omega^{-lk})t_{j}^{\lbr a^k\rbr}+(-1)^{l}t_{j}^{\lbr a^{n-2}\rbr}
\eea
for $0\leq l\leq n-2$, together with $t_{j}^{\lbr b\rbr}$ and $t_{j}^{\lbr ab\rbr}$. The coordinates transforming matrix is obvious.\\

Now we rewrite the operator $\big(\frac{A_{p+1}(V_{\rho_{1}})z^p}{(p+1)!}\big)^{\wedge}$ in the new coordinates as
\ben
&&\big(\frac{A_{p+1}(V_{\rho_{1}})z^p}{(p+1)!}\big)^{\wedge}\\
&=&\frac{2B_{p+1}}{(p+1)!}\pd_{\lbr 1\rbr,1+p}
-\frac{2B_{p+1}}{(p+1)!}\sum_{l= 0}^{\infty}t_{l}^{\lbr 1\rbr}\pd_{\lbr 1\rbr,l+p}-\sum_{k=1}^{n-3}\frac{B_{p+1}(\frac{k}{2n-4})+B_{p+1}(\frac{2n-4-k}{2n-4})}{(p+1)!}\sum_{l= 0}^{\infty}t_{l}^{\lbr a^k\rbr}\pd_{\lbr a^k\rbr,l+p}\\
&&-\frac{2B_{p+1}(\frac{n-2}{2n-4})}{(p+1)!}\sum_{l=0}^{\infty}t_{l}^{\lbr a^{n-2}\rbr}\pd_{\lbr a^{n-2}\rbr,l+p}-\frac{B_{p+1}(\frac{1}{4})+B_{p+1}(\frac{3}{4})}{(p+1)!}\sum_{l=0}^{\infty}(t_{l}^{\lbr b\rbr}\pd_{\lbr b\rbr,l+p}+t_{l}^{\lbr ab\rbr}\pd_{\lbr ab\rbr,l+p})
\\
&&+\frac{\hbar ^2}{2}\sum_{l= 0}^{p-1}(-1)^l\Bigg((4n-8)\frac{2B_{p+1}}{(p+1)!}\pd_{\lbr 1\rbr,l}\pd_{\lbr 1\rbr,p-1-l}+(2n-4)\sum_{k=1}^{n-3}\frac{B_{p+1}(\frac{k}{2n-4})+B_{p+1}(\frac{2n-4-k}{2n-4})}{(p+1)!}\\
&& \cdot\pd_{\lbr a^k \rbr,l}\pd_{\lbr a^k\rbr,p-1-l}+(4n-8)\frac{2B_{p+1}(\frac{1}{2})}{(p+1)!}\pd_{\lbr a^{n-2}\rbr,l}\pd_{\lbr a^{n-2}\rbr,p-1-l}+
4\cdot\frac{B_{p+1}(\frac{1}{4})+B_{p+1}(\frac{3}{4})}{(p+1)!}\\
&&\cdot(\pd_{\lbr b\rbr,l}\pd_{
\lbr b \rbr,p-1-l}+\pd_{\lbr ab\rbr,l}\pd_{
\lbr ab \rbr,p-1-l})\Bigg)\\
&=&
\frac{2B_{p+1}}{(p+1)!}\Big(\sum_{i=1}^{4}\frac{\pd}{\pd u_{1+p}^{\psi_{i}}}+\sum_{j=1}^{n-3}\frac{\pd}{\pd u_{1+p}^{\rho_{j}}}\Big)
-\frac{2B_{p+1}}{(p+1)!}\frac{1}{4n-8}\sum_{l= 0}^{\infty}\Big(\sum_{i=1}^{4}u_{l}^{\psi_i}+4\sum_{j=1}^{n-3}u_{l}^{\rho_{j}}\Big)\\
&&\Big(\sum_{i=1}^{4}\frac{\pd}{\pd u_{l+p}^{\psi_{i}}}+\sum_{j=1}^{n-3}\frac{\pd}{\pd u_{l+p}^{\rho_{j}}}\Big)-\sum_{k=1}^{n-3}\frac{B_{p+1}(\frac{k}{2n-4})+B_{p+1}(\frac{2n-4-k}{2n-4})}{(p+1)!}\frac{1}{4n-8}\\
&&\cdot\sum_{l= 0}^{\infty}\Big( u_{l}^{\psi_{1}}+u_{l}^{\psi_{2}}+(-1)^{k}u_{l}^{\psi_{3}}+(-1)^{k}u_{l}^{\psi_{4}}+
2\sum_{j=1}^{n-3}(\omega^{kj}+\omega^{-kj})u_{l}^{\rho_{j}}\Big) \\
&& \cdot\Big(2\frac{\pd}{\pd u_{p+l}^{\psi_{1}}}+2\frac{\pd}{\pd u_{p+l}^{\psi_{2}}}+2(-1)^{k}\frac{\pd}{\pd u_{p+l}^{\psi_{3}}}+2(-1)^{k}\frac{\pd}{\pd u_{p+l}^{\psi_{4}}}+\sum_{j=1}^{n-3}(\omega^{kj}+\omega^{-kj})\frac{\pd}{\pd u_{p+l}^{\rho_j}}\Big)\\
&&-\frac{2B_{p+1}(\frac{1}{2})}{(p+1)!}
\frac{1}{4n-8}\sum_{l=0}^{\infty}\Big(\sum_{i=1}^{4}u_{l}^{\psi_i}
+4\sum_{j=1}^{n-3}(-1)^{j}u_{l}^{\rho_{j}}\Big)\Big(\sum_{i=1}^{4}\frac{\pd}{\pd u_{l+p}^{\psi_{i}}}+\sum_{j=1}^{n-3}(-1)^{j}\frac{\pd}{\pd u_{l+p}^{\rho_{j}}}\Big)\\
&&-\frac{B_{p+1}(\frac{1}{4})+B_{p+1}(\frac{3}{4})}{(p+1)!}\sum_{l=0}^{\infty}\Bigg( \frac{1}{4}(u_{l}^{\psi_{1}} -u_{l}^{\psi_{2}}-u_{l}^{\psi_{3}} +u_{l}^{\psi_{4}})(\frac{\pd}{\pd u_{p+l}^{\psi_1}}-\frac{\pd}{\pd u_{p+l}^{\psi_2}}-\frac{\pd}{\pd u_{p+l}^{\psi_3}}+\frac{\pd}{\pd u_{p+l}^{\psi_4}})\\
&&+\frac{1}{4}(u_{l}^{\psi_{1}} -u_{l}^{\psi_{2}}+u_{l}^{\psi_{3}} -u_{l}^{\psi_{4}})(\frac{\pd}{\pd u_{p+l}^{\psi_1}}-\frac{\pd}{\pd u_{p+l}^{\psi_2}}+\frac{\pd}{\pd u_{p+l}^{\psi_3}}-\frac{\pd}{\pd u_{p+l}^{\psi_4}})\Bigg)
\\
&&+\frac{\hbar ^2}{2}\sum_{l= 0}^{p-1}(-1)^l\Bigg[(4n-8)\frac{2B_{p+1}}{(p+1)!}\Big(\sum_{i=1}^{4}\frac{\pd}{\pd u_{l}^{\psi_{i}}}+\sum_{j=1}^{n-3}\frac{\pd}{\pd u_{l}^{\rho_{j}}}\Big)\Big(\sum_{i=1}^{4}\frac{\pd}{\pd u_{p-1-l}^{\psi_{i}}}+\sum_{j=1}^{n-3}\frac{\pd}{\pd u_{p-1-l}^{\rho_{j}}}\Big)\\
&&+(2n-4)\sum_{k=1}^{n-3}\frac{B_{p+1}(\frac{k}{2n-4})+B_{p+1}(\frac{2n-4-k}{2n-4})}{(p+1)!}\\
&&\cdot 4\Big(\frac{\pd}{\pd u_{l}^{\psi_{1}}}+\frac{\pd}{\pd u_{l}^{\psi_{2}}}+(-1)^{k}\frac{\pd}{\pd u_{l}^{\psi_{3}}}+(-1)^{k}\frac{\pd}{\pd u_{l}^{\psi_{4}}}+\sum_{j=1}^{n-3}\frac{\omega^{kj}+\omega^{-kj}}{2}\frac{\pd}{\pd u_{l}^{\rho_j}}\Big)\\
&&\cdot\Big(\frac{\pd}{\pd u_{p-1-l}^{\psi_{1}}}+\frac{\pd}{\pd u_{p-1-l}^{\psi_{2}}}+(-1)^{k}\frac{\pd}{\pd u_{p-1-l}^{\psi_{3}}}+(-1)^{k}\frac{\pd}{\pd u_{p-1-l}^{\psi_{4}}}+\sum_{j=1}^{n-3}\frac{\omega^{kj}+\omega^{-kj}}{2}\frac{\pd}{\pd u_{p-1-l}^{\rho_j}}\Big)\\
&&+(4n-8)\cdot\frac{2B_{p+1}(\frac{1}{2})}{(p+1)!}\Big(\sum_{i=1}^{4}\frac{\pd}{\pd u_{l}^{\psi_{i}}}+\sum_{j=1}^{n-3}(-1)^{j}\frac{\pd}{\pd u_{l}^{\rho_{j}}}\Big)\Big(\sum_{i=1}^{4}\frac{\pd}{\pd u_{p-1-l}^{\psi_{i}}}+\sum_{j=1}^{n-3}(-1)^{j}\frac{\pd}{\pd u_{p-1-l}^{\rho_{j}}}\Big)\\
&&+4(n-2)^{2}\cdot\frac{B_{p+1}(\frac{1}{4})+B_{p+1}(\frac{3}{4})}{(p+1)!}\Big((\frac{\pd}{\pd u_{l}^{\psi_1}}-\frac{\pd}{\pd u_{l}^{\psi_2}}-\frac{\pd}{\pd u_{l}^{\psi_3}}+\frac{\pd}{\pd u_{l}^{\psi_4}})\\
&&\cdot (\frac{\pd}{\pd u_{p-1-l}^{\psi_1}}-\frac{\pd}{\pd u_{p-1-l}^{\psi_2}}-\frac{\pd}{\pd u_{p-1-l}^{\psi_3}}+\frac{\pd}{\pd u_{p-1-l}^{\psi_4}})\\
&&+(\frac{\pd}{\pd u_{l}^{\psi_1}}-\frac{\pd}{\pd u_{l}^{\psi_2}}+\frac{\pd}{\pd u_{l}^{\psi_3}}-\frac{\pd}{\pd u_{l}^{\psi_4}})(\frac{\pd}{\pd u_{p-1-l}^{\psi_1}}-\frac{\pd}{\pd u_{p-1-l}^{\psi_2}}+\frac{\pd}{\pd u_{p-1-l}^{\psi_3}}-\frac{\pd}{\pd u_{p-1-l}^{\psi_4}})\Big)\Bigg]\\
&=&
\frac{2B_{p+1}}{(p+1)!}\Big(\sum_{i=1}^{4}\frac{\pd}{\pd u_{1+p}^{\psi_{i}}}+\sum_{j=1}^{n-3}\frac{\pd}{\pd u_{1+p}^{\rho_{j}}}\Big)
\\&&-\frac{1}{4n-8}\sum_{l=0}^{\infty}\Bigg[\frac{2\beta_{p+1}(2n-4,0)}{(p+1)!}\Big((u_{l}^{\psi_{1}}+u_{l}^{\psi_{2}})(\frac{\pd}{\pd u_{l+p}^{\psi_{1}}}+\frac{\pd}{\pd u_{l+p}^{\psi_{2}}})
+(u_{l}^{\psi_{3}}+u_{l}^{\psi_{4}})(\frac{\pd}{\pd u_{l+p}^{\psi_{3}}}+\frac{\pd}{\pd u_{l+p}^{\psi_{4}}})\Big)\\
&&
+\frac{2\beta_{p+1}(2n-4,n-2)}{(p+1)!}\Big((u_{l}^{\psi_{1}}+u_{l}^{\psi_{2}})(\frac{\pd}{\pd u_{l+p}^{\psi_{3}}}+\frac{\pd}{\pd u_{l+p}^{\psi_{4}}})+(u_{l}^{\psi_{3}}+u_{l}^{\psi_{4}})(\frac{\pd}{\pd u_{l+p}^{\psi_{1}}}+\frac{\pd}{\pd u_{l+p}^{\psi_{2}}})\Big)\\
&&+\sum_{j=1}^{n-3}\frac{\beta_{p+1}(2n-4,j)+\beta_{p+1}(2n-4,-j)}{(p+1)!}(u_{l}^{\psi_{1}}+u_{l}^{\psi_{2}})\frac{\pd}{\pd u_{l+p}^{\rho_{j}}}\\
&&+\sum_{j=1}^{n-3}\frac{\beta_{p+1}(2n-4,n-2+j)+\beta_{p+1}(2n-4,n-2-j)}{(p+1)!}(u_{l}^{\psi_{3}}+u_{l}^{\psi_{4}})\frac{\pd}{\pd u_{l+p}^{\rho_{j}}}\\
&&+4\sum_{j=1}^{n-3}\frac{\beta_{p+1}(2n-4,j)+\beta_{p+1}(2n-4,-j)}{(p+1)!}u_{l}^{\rho_{j}}(\frac{\pd}{\pd u_{l+p}^{\psi_{1}}}+\frac{\pd}{\pd u_{l+p}^{\psi_{2}}})\\
&&+4\sum_{j=1}^{n-3}\frac{\beta_{p+1}(2n-4,n-2+j)+\beta_{p+1}(2n-4,n-2-j)}{(p+1)!}u_{l}^{\rho_{j}}(\frac{\pd}{\pd u_{l+p}^{\psi_{3}}}+\frac{\pd}{\pd u_{l+p}^{\psi_{4}}})\\
&&+2\sum_{i=1}^{n-3}\sum_{j=1}^{n-3}\Big(\frac{\beta_{p+1}(2n-4,i+j)+\beta_{p+1}(2n-4,i-j)}{(p+1)!}\\
&&+\frac{\beta_{p+1}(2n-4,-i+j)+\beta_{p+1}(2n-4,-i-j)}{(p+1)!}\Big)
u_{l}^{\rho_{i}}\frac{\pd}{\pd u_{l+p}^{\rho_{j}}}
\Bigg]\\
&&-\sum_{l=0}^{\infty}\frac{B_{p+1}(\frac{1}{4})+B_{p+1}(\frac{3}{4})}{(p+1)!}\Bigg( \frac{1}{2}(u_{l}^{\psi_{1}} -u_{l}^{\psi_{2}})(\frac{\pd}{\pd u_{p+l}^{\psi_1}}-\frac{\pd}{\pd u_{p+l}^{\psi_2}})+\frac{1}{2}(u_{l}^{\psi_{3}} -u_{l}^{\psi_{4}})(\frac{\pd}{\pd u_{p+l}^{\psi_3}}-\frac{\pd}{\pd u_{p+l}^{\psi_4}})\Bigg)\\
&&+\frac{\hbar^2}{2}\sum_{l=0}^{p-1}(-1)^{l}\Bigg[8(n-2)\frac{\beta_{p+1}(2n-4,0)}{(p+1)!}
\Big( (\frac{\pd}{\pd u_{l}^{\psi_1}}+\frac{\pd}{\pd u_{l}^{\psi_2}})(\frac{\pd}{\pd u_{p-1-l}^{\psi_1}}+\frac{\pd}{\pd u_{p-1-l}^{\psi_2}})\\
&&+(\frac{\pd}{\pd u_{l}^{\psi_3}}+\frac{\pd}{\pd u_{l}^{\psi_4}})(\frac{\pd}{\pd u_{p-1-l}^{\psi_3}}+\frac{\pd}{\pd u_{p-1-l}^{\psi_4}})\Big)+8(n-2)\frac{\beta_{p+1}(2n-4,n-2)}{(p+1)!}\\
&&\cdot\Big( (\frac{\pd}{\pd u_{l}^{\psi_1}}+\frac{\pd}{\pd u_{l}^{\psi_2}})(\frac{\pd}{\pd u_{p-1-l}^{\psi_3}}+\frac{\pd}{\pd u_{p-1-l}^{\psi_4}})+(\frac{\pd}{\pd u_{l}^{\psi_3}}+\frac{\pd}{\pd u_{l}^{\psi_4}})(\frac{\pd}{\pd u_{p-1-l}^{\psi_1}}+\frac{\pd}{\pd u_{p-1-l}^{\psi_2}})\Big )\\
&&+8(n-2)\sum_{j=1}^{n-3}\frac{\beta_{p+1}(2n-4,j)+\beta_{p+1}(2n-4,-j)}{(p+1)!}
(\frac{\pd}{\pd u_{l}^{\psi_1}}+\frac{\pd}{\pd u_{l}^{\psi_2}})\frac{\pd}{\pd u_{p-1-l}^{\rho_j}}\\
&&+8(n-2)\sum_{j=1}^{n-3}\frac{\beta_{p+1}(2n-4,n-2+j)+\beta_{p+1}(2n-4,n-2-j)}{(p+1)!}
(\frac{\pd}{\pd u_{l}^{\psi_3}}+\frac{\pd}{\pd u_{l}^{\psi_4}})\frac{\pd}{\pd u_{p-1-l}^{\rho_j}}\\
&&+(2n-4)\frac{1}{(p+1)!}\sum_{i=1}^{n-3}\sum_{j=1}^{n-3}\Big(\beta_{p+1}(2n-4,i+j)+\beta_{p+1}(2n-4,i-j)\\
&&+\beta_{p+1}(2n-4,-i+j)+\beta_{p+1}(2n-4,-i-j)\Big )\frac{\pd}{\pd u_{l}^{\rho_i}}\frac{\pd}{\pd u_{p-1-l}^{\rho_j}}\\
&&+8(n-2)^{2}\cdot\frac{B_{p+1}(\frac{1}{4})+B_{p+1}(\frac{3}{4})}{(p+1)!}\Bigg((\frac{\pd}{\pd u_{l}^{\psi_{1}}} -\frac{\pd}{\pd u_{l}^{\psi_{2}}})(\frac{\pd}{\pd u_{p-1-l}^{\psi_1}}-\frac{\pd}{\pd u_{p-1-l}^{\psi_2}})\\
&&+(\frac{\pd}{\pd u_{l}^{\psi_{3}}}-\frac{\pd}{\pd u_{l}^{\psi_{4}}})(\frac{\pd}{\pd u_{p-1-l}^{\psi_3}}-\frac{\pd}{\pd u_{p-1-l}^{\psi_4}})\Bigg)
\Bigg]\\
&=&\frac{2B_{p+1}}{(p+1)!}\sum_{j=0}^{n-2}\frac{\pd}{\pd u_{1+p}^{\rho_{j}}}
\\&&-\frac{1}{4n-8}\sum_{l=0}^{\infty}\Bigg[\frac{4\beta_{p+1}(2n-4,0)}{(p+1)!}\Big(u_{l}^{\rho_{0}}\frac{\pd}{\pd u_{l+p}^{\rho_{0}}}
+u_{l}^{\rho_{n-2}}\frac{\pd}{\pd u_{l+p}^{\rho_{n-2}}}\Big)\\
&&
+\frac{4\beta_{p+1}(2n-4,n-2)}{(p+1)!}\Big(u_{l}^{\rho_{0}}\frac{\pd}{\pd u_{l+p}^{\rho_{n-2}}}
+u_{l}^{\rho_{n-2}}\frac{\pd}{\pd u_{l+p}^{\rho_{0}}}\Big)\\
&&+2\sum_{j=1}^{n-3}\frac{\beta_{p+1}(2n-4,j)+\beta_{p+1}(2n-4,-j)}{(p+1)!}u_{l}^{\rho_{0}}\frac{\pd}{\pd u_{l+p}^{\rho_{j}}}\\
&&+2\sum_{j=1}^{n-3}\frac{\beta_{p+1}(2n-4,n-2+j)+\beta_{p+1}(2n-4,n-2-j)}{(p+1)!}u_{l}^{\rho_{n-2}}\frac{\pd}{\pd u_{l+p}^{\rho_{j}}}\\
&&+4\sum_{j=1}^{n-3}\frac{\beta_{p+1}(2n-4,j)+\beta_{p+1}(2n-4,-j)}{(p+1)!}u_{l}^{\rho_{j}}\frac{\pd}{\pd u_{l+p}^{\rho_{0}}}\\
&&+4\sum_{j=1}^{n-3}\frac{\beta_{p+1}(2n-4,n-2+j)+\beta_{p+1}(2n-4,n-2-j)}{(p+1)!}u_{l}^{\rho_{j}}\frac{\pd}{\pd u_{l+p}^{\rho_{n-2}}}\\
&&+2\sum_{i}^{n-3}\sum_{i}^{n-3}\Big(\frac{\beta_{p+1}(2n-4,i+j)+\beta_{p+1}(2n-4,i-j)}{(p+1)!}\\
&&+\frac{\beta_{p+1}(2n-4,-i+j)+\beta_{p+1}(2n-4,-i-j)}{(p+1)!}\Big)
u_{l}^{\rho_{i}}\frac{\pd}{\pd u_{l+p}^{\rho_{j}}}
\Bigg]\\
&&-\frac{B_{p+1}(\frac{1}{4})+B_{p+1}(\frac{3}{4})}{(p+1)!}\sum_{l=0}^{\infty}(t_{l}^{\lbr b\rbr}\pd_{\lbr b\rbr,l+p}+t_{l}^{\lbr ab\rbr}\pd_{\lbr ab\rbr,l+p})\\
&&+\frac{\hbar^2}{2}\sum_{l=0}^{p-1}(-1)^{l}\Bigg[(8n-2)\frac{\beta_{p+1}(2n-4,0)}{(p+1)!}
\Big( \frac{\pd}{\pd u_{l}^{\rho_0}}\frac{\pd}{\pd u_{p-1-l}^{\rho_0}}+\frac{\pd}{\pd u_{l}^{\rho_{n-2}}}\frac{\pd}{\pd u_{p-1-l}^{\rho_{n-2}}}\Big)\\
&&+(8n-2)\frac{\beta_{p+1}(2n-4,n-2)}{(p+1)!}\Big( \frac{\pd}{\pd u_{l}^{\rho_{0}}}\frac{\pd}{\pd u_{p-1-l}^{\rho_{n-2}}}+\frac{\pd}{\pd u_{l}^{\rho_{n-2}}}\frac{\pd}{\pd u_{p-1-l}^{\rho_0}}\Big )\\
&&+(8n-2)\sum_{j=1}^{n-3}\frac{\beta_{p+1}(2n-4,j)+\beta_{p+1}(2n-4,-j)}{(p+1)!}
\frac{\pd}{\pd u_{l}^{\rho_{0}}}\frac{\pd}{\pd u_{p-1-l}^{\rho_j}}\\
&&+(8n-2)\sum_{j=1}^{n-3}\frac{\beta_{p+1}(2n-4,n-2+j)+\beta_{p+1}(2n-4,n-2-j)}{(p+1)!}
\frac{\pd}{\pd u_{l}^{\rho_{n-2}}}\frac{\pd}{\pd u_{p-1-l}^{\rho_j}}\\
&&+(2n-4)\frac{1}{(p+1)!}\sum_{i=1}^{n-3}\sum_{j=1}^{n-3}\Big(\beta_{p+1}(2n-4,i+j)+\beta_{p+1}(2n-4,i-j)\\
&&+\beta_{p+1}(2n-4,-i+j)+\beta_{p+1}(2n-4,-i-j)\Big )\frac{\pd}{\pd u_{l}^{\rho_i}}\frac{\pd}{\pd u_{p-1-l}^{\rho_j}}\\
&&+4\cdot\frac{B_{p+1}(\frac{1}{4})+B_{p+1}(\frac{3}{4})}{(p+1)!}(\pd_{\lbr b\rbr,l}\pd_{
\lbr b \rbr,p-1-l}+\pd_{\lbr ab\rbr,l}\pd_{
\lbr ab \rbr,p-1-l})
\Bigg].\\
\een

Note that when $2\nmid (p+1)$ and $p\geq 2$, $\big(\frac{A_{p+1}(V_{\rho_{1}})z^p}{(p+1)!}\big)^{\wedge}=0$. In the following we assume $2|(p+1)$. Let $u_{l}^{\rho_{j}}=\frac{1}{2}(v_{l}^{j}+v_{l}^{2n-4-j})$  and $\bar{u}_{l}^{\rho_{j}}=\frac{1}{2}(v_{l}^{j}-v_{l}^{2n-4-j})$, for $1\leq j\leq n-3$. We \emph{formally modify} the operator $\big(\frac{A_{p+1}(V_{\rho_{1}})z^p}{(p+1)!}\big)^{\wedge}$ as
\bea\label{109}
&&\big(\frac{A_{p+1}(V_{\rho_{1}})z^p}{(p+1)!}\big)_{\text{mod}}^{\wedge}\nonumber\\
&=&\frac{2B_{p+1}}{(p+1)!}\sum_{j=0}^{2n-5}\frac{\pd}{\pd v_{1+p}^{j}}-\frac{4}{4n-8}\sum_{l=0}^{\infty}
\sum_{i=0}^{2n-5}\sum_{j=0}^{2n-5}\frac{\beta_{p+1}(2n-4,j-i)}{(p+1)!}v_{l}^{i}\frac{\pd}{\pd v_{l+p}^{j}}\nonumber\\
&&-\frac{B_{p+1}(\frac{1}{4})+B_{p+1}(\frac{3}{4})}{(p+1)!}\sum_{l=0}^{\infty}(t_{l}^{\lbr b\rbr}\pd_{\lbr b\rbr,l+p}+t_{l}^{\lbr ab\rbr}\pd_{\lbr ab\rbr,l+p})\nonumber\\
&&+\frac{\hbar^2}{2}\sum_{l=0}^{p-1}(-1)^{l}\Bigg[8(n-2)\sum_{i=0}^{2n-5}\sum_{j=0}^{2n-5}\frac{\beta_{p+1}(2n-4,i+j)}{(p+1)!}\frac{\pd}{\pd v_{l}^{i}}\frac{\pd}{\pd v_{p-1-l}^{j}}\nonumber\\
&&+4\cdot\frac{B_{p+1}(\frac{1}{4})+B_{p+1}(\frac{3}{4})}{(p+1)!}(\pd_{\lbr b\rbr,l}\pd_{
\lbr b \rbr,p-1-l}+\pd_{\lbr ab\rbr,l}\pd_{
\lbr ab \rbr,p-1-l})
\Bigg].\nonumber\\
\eea
By \cite{JK}, the all genera free energy of $B\hat{D}_n$
\bea\label{56}
F^{\hat{D}_n}&=&F\Bigg(\Big\{\frac{v_{l}^{0}+(n-2)t_{l}^{\lbr b\rbr}+(n-2)t_{l}^{\lbr ab\rbr}
}{(4n-8)^{\frac{2(1-l)}{3}}}\Big\}\Bigg)+F\Bigg(\Big\{\frac{v_{l}^{0}-(n-2)t_{l}^{\lbr b\rbr}-(n-2)t_{l}^{\lbr ab\rbr}
}{(4n-8)^{\frac{2(1-l)}{3}}}\Big\}\Bigg)\nonumber\\
&&+ F\Bigg(\Big\{\frac{v_{l}^{n-2}-(n-2)t_{l}^{\lbr b\rbr}+(n-2)t_{l}^{\lbr ab\rbr}
}{(4n-8)^{\frac{2(1-l)}{3}}}\Big\}\Bigg) +F\Bigg(\Big\{\frac{v_{l}^{n-2}+(n-2)t_{l}^{\lbr b\rbr}-(n-2)t_{l}^{\lbr ab\rbr}
}{(4n-8)^{\frac{2(1-l)}{3}}}\Big\}\Bigg)\nonumber\\
&&+\sum_{j=1}^{n-3}F\Bigg(\Big\{\frac{v_{l}^{j}+v_{l}^{2n-4-j}}{(2n-4)^{
\frac{2(1-l)}{3}}}\Big\}\Bigg).
\eea
Recall the quantum Riemann-Roch theorem and the notations in section 1. We have
\begin{lemma}\label{35}
\ben
\big(\frac{A_{p+1}(V_{\rho_{1}})z^p}{(p+1)!}\big)^{\wedge}Z^{\mathbf{c},\rho}=
\big(\frac{A_{p+1}(V_{\rho_{1}})z^p}{(p+1)!}\big)_{\emph{mod}}^{\wedge}Z^{\mathbf{c},\rho}.
\een
\end{lemma}
Proof:  Under the change of coordinates $\{u_{l}^{\rho_{j}}, \bar{u}_{l}^{\rho_{j}}\}$ to $\{v_{l}^{j}, v_{l}^{2n-4-j}\}$ we have
\ben
(a+b)\frac{\pd}{\pd u_{l}^{\rho_{j}}}&=&(a+b)\Big(\frac{\pd}{\pd v_{l}^{j}}+\frac{\pd}{\pd v_{l}^{2n-4-j}}\Big)\\
&=& 2a\frac{\pd}{\pd v_{l}^{j}}+ 2b\frac{\pd}{\pd v_{l}^{2n-4-j}}-(a-b)\Big(\frac{\pd}{\pd v_{l}^{j}}-\frac{\pd}{\pd v_{l}^{2n-4-j}}\Big).
\een
Thus when  $f$ is a function of $u_{l}^{\rho_{j}}$,
\ben
(a+b)\frac{\pd}{\pd u_{l}^{\rho_{j}}}f= \Big(2a\frac{\pd}{\pd v_{l}^{j}}+ 2b\frac{\pd}{\pd v_{l}^{2n-4-j}}\Big)f.
\een
This elementary fact applies in our case. \hfill\qedsymbol

\begin{theorem}\label{70}
For fixed $m\geq 2$, $\langle e_{\lbr b \rbr}^{2m} \rangle^{[\mathbb{C}^2 /\hat{D}_n]}$ is a Laurent polynomial of $n-2$.
\end{theorem}
Proof: We separate the operator $\big(\frac{A_{p+1}(V_{\rho_{1}})z^p}{(p+1)!}\big)_{\text{mod}}^{\wedge}$ into two parts,
\ben
L_{p}&:=&\frac{2B_{p+1}}{(p+1)!}\sum_{j=0}^{2n-5}\frac{\pd}{\pd v_{1+p}^{j}}-\frac{1}{n-2}\sum_{l=0}^{\infty}
\sum_{i=0}^{2n-5}\sum_{j=0}^{2n-5}\frac{\beta_{p+1}(2n-4,j-i)}{(p+1)!}v_{l}^{i}\frac{\pd}{\pd v_{l+p}^{j}}\\
&&+\frac{\hbar^2}{2}\sum_{l=0}^{p-1}(-1)^{l}8(n-2)\sum_{i=0}^{2n-5}\sum_{j=0}^{2n-5}\frac{\beta_{p+1}(2n-4,i+j)}{(p+1)!}\frac{\pd}{\pd v_{l}^{i}}\frac{\pd}{\pd v_{p-1-l}^{j}},\\
M_{p}&:=&-\frac{B_{p+1}(\frac{1}{4})+B_{p+1}(\frac{3}{4})}{(p+1)!}\sum_{l=0}^{\infty}(t_{l}^{\lbr b\rbr}\pd_{\lbr b\rbr,l+p}+t_{l}^{\lbr ab\rbr}\pd_{\lbr ab\rbr,l+p})\\
&&+\frac{\hbar^2}{2}\sum_{l=0}^{p-1}(-1)^{l}4\cdot\frac{B_{p+1}(\frac{1}{4})+B_{p+1}(\frac{3}{4})}{(p+1)!}(\pd_{\lbr b\rbr,l}\pd_{
\lbr b \rbr,p-1-l}+\pd_{\lbr ab\rbr,l}\pd_{
\lbr ab \rbr,p-1-l}).\\
\een
Thus $L_{p_1}$ commutes with $M_{p_2}$ for  $p_1, p_2\geq 1$. By theorem \ref{18}, it suffices to show that, for any $p_1,\cdots, p_k\geq 1$,
the coefficient of every term $t_{l_1}^{\lbr b\rbr}\cdots t_{l_{r_1}}^{\lbr b\rbr}t_{l_{r_1+1}}^{\lbr ab\rbr}\cdots t_{l_{r_2+1}}^{\lbr ab\rbr}$ in the  \emph{LTHE contributions} in
\ben
\Big(L_{p_1}\cdots L_{p_k}F_{0}^{\hat{D}_n}\Big)\Big|_{\text{all}\hspace{0.2cm} v_{l}^{j}=0}
\een
 is a Laurent polynomial of $n-2$. But we have
 \ben
&& \frac{2B_{k_{1}+1}}{(k_{1}+1)!}\frac{1}{n-2}\sum_{j_{1}=0}^{2n-5}\sum_{l=0}^{\infty}
\sum_{i_2=0}^{2n-5}\sum_{j_2=0}^{2n-5}\frac{\beta_{k_{2}+1}(2n-4,j_{2}-i_{2})}{(k_2+1)!}\Big(\frac{\pd}{\pd v_{1+k_{1}}^{j_1}}v_{l}^{i_{2}}\Big)\frac{\pd}{\pd v_{l+k_2}^{j_2}}\\
&=& \frac{2B_{k_{1}+1}}{(k_{1}+1)!}\frac{1}{n-2}\sum_{j_{1}=0}^{2n-5}
\sum_{j_2=0}^{2n-5}\frac{\beta_{k_{2}+1}(2n-4,j_{2}-j_{1})}{(k_2+1)!}\frac{\pd}{\pd v_{k_1+k_2+1}^{j_2}}\\
&=& \frac{2B_{k_{1}+1}}{(k_{1}+1)!}\frac{1}{n-2}
\sum_{j_2=0}^{2n-5}\frac{(2n-4)B_{k_2+1}}{(k_2+1)!}\frac{\pd}{\pd v_{k_1+k_2+1}^{j_2}}\\
&=& \frac{4B_{k_{1}+1}B_{k_2+1}}{(k_{1}+1)!(k_2+1)!}
\sum_{j=0}^{2n-5}\frac{\pd}{\pd v_{k_1+k_2+1}^{j}},
 \een
and
\ben
&&\sum_{l_1=0}^{k_1-1}(-1)^{l_1}8(n-2)\sum_{i_1=0}^{2n-5}\sum_{j_1=0}^{2n-5}\frac{1}{n-2}\sum_{l_2=0}^{\infty}
\sum_{i_2=0}^{2n-5}\sum_{j_2=0}^{2n-5}\frac{\beta_{k_1+1}(2n-4,i_1+j_1)}{(k_1+1)!}\\
&&\cdot\frac{\beta_{k_2+1}(2n-4,j_2-i_2)}{(k_2+1)!}\frac{\pd}{\pd v_{l_1}^{i_1}}\Big(\frac{\pd}{\pd v_{k_1-1-l_1}^{j_1}}v_{l_2}^{i_2}\Big)\frac{\pd}{\pd v_{l_2+k_2}^{j_2}}\\
&=&\sum_{l_1=0}^{k_1-1}(-1)^{l_1}8\sum_{i_1=0}^{2n-5}\sum_{j_1=0}^{2n-5}
\sum_{j_2=0}^{2n-5}\frac{\beta_{k_1+1}(2n-4,i_1+j_1)}{(k_1+1)!}\frac{\beta_{k_2+1}(2n-4,j_2-j_1)}{(k_2+1)!}\frac{\pd}{\pd v_{l_1}^{i_1}}\frac{\pd}{\pd v_{k_1+k_2-1-l_1}^{j_2}}\\
&=&\sum_{l=0}^{k_1-1}(-1)^{l}16(n-2)\sum_{i=0}^{2n-5}
\sum_{j=0}^{2n-5}\frac{\beta_{k_1+1,k_2+1}(2n-4,i+j)}{(k_1+1)!(k_2+1)!}\frac{\pd}{\pd v_{l}^{i}}\frac{\pd}{\pd v_{k_1+k_2-1-l}^{j}}.
\een

 \ben
 &&\sum_{j=0}^{n-1}\beta_{M}(n,i+j)\beta_{N}(k-j)\\
 &=& \sum_{j=0}^{n-1}\sum_{a=0}^{n-1}\sum_{b=0}^{n-1}\zeta_{n}^{a(i+j)}\zeta_{n}^{b(k-j)}B_{M}(\frac{a}{n})B_{N}(\frac{b}{n})\\
 &=& n\sum_{a=0}^{n-1}\zeta_{n}^{a(i+k)}B_{M}(\frac{a}{n})B_{N}(\frac{a}{n})\\
 &=& n\beta_{M\coprod N}(n,i+k).
 \een
Since we shall take all $v_{l}^{j}=0$ at last, which forces that every $v_{l}^{j}$ must be differentiated. Furthermore, We have the \\

\textbf{Claim 1}: it suffices consider the operators of the form
\bea\label{40}
\sum_{j=0}^{2n-5}\frac{\pd}{\pd v_{1+p}^{j}}
\eea
or
\bea\label{41}
\frac{\hbar^2}{2}\sum_{i=0}^{2n-5}
\sum_{j=0}^{2n-5}\beta_{M}(2n-4,i+j)\frac{\pd}{\pd v_{l}^{i}}\frac{\pd}{\pd v_{k}^{j}}.
\eea
Proof of the claim 1: It is obvious, by noting that the ranges of the subscripts of $v_{l}^{i}$ are independent of $n$.\hfill\qedsymbol\\

Now let us consider only the operators of the form (\ref{41}). In this case every LTHE concerned has no half edges. Because of the existence of the set of subscripts $M$, every LTHE naturally equipped with the $M$'s becomes a \emph{decorated Feymann diagram}. Fix such a decorated Feymann diagram $\Gamma$, we shall show that \\

\textbf{Claim 2}: $\text{Cont}(\Gamma)$ can be written in the form $\sum_{i}f_i S_{n-2}(\Gamma_{i})$, where $f_i$ are Laurent polynomials of $n-2$, and $\Gamma_{i}$  run over the set of subtree of $\Gamma$.\\

Proof of Claim 2: We prove this by induction along the stratification of $\Gamma$. First let $v_1$ be a vertex in $V(\Gamma(1))\backslash V(\Gamma(2))$. Suppose $\Gamma$ is not a star diagram. Then the vertices neighboring to $v_1$ are all leaves, except exactly one. Let the vertices neighboring to $v_1$ be $u_1,\cdots, u_r, u_{r+1}$, the corresponding edges be $e_1,\cdots, e_r, e_{r+1}$, and assume that $u_1,\cdots,u_{r}$ are leaves, and thus $u_{r+1}$ is not. Let $M_{i}$ be the decoration of $e_i$, $1\leq i\leq r+1$.
For every leaf neighboring to $v_1$, the first derivative $(\frac{\pd}{\pd v_{l}^{i}}F_{0}^{\hat{D}_{n}})\big|_{\text{all}\hspace{0.1cm} v_{k}^{j}=0}$ gives no contribution, unless $i=0$ or $n-2$. Thus the contribution of such a leaf is
\bea\label{48}
\beta_{M_i}(2n-4,j)+\beta_{M_i}(2n-4,n-2+j)
\eea
multiplied by
\bea\label{49}
(\pd_{l_i}F_{0}^{\hat{D}_{n}})\Big(\Big\{\frac{(n-2)t_{l}^{\lbr b\rbr}+(n-2)t_{l}^{\lbr ab\rbr}
}{(4n-8)^{\frac{2(1-l)}{3}}}\Big\}\Big)+(\pd_{l_i}F_{0}^{\hat{D}_{n}})\Big(\Big\{\frac{-(n-2)t_{l}^{\lbr b\rbr}-(n-2)t_{l}^{\lbr ab\rbr}
}{(4n-8)^{\frac{2(1-l)}{3}}}\Big\}\Big).
\eea
But the factor (\ref{49}) has all coefficients Laurent polynomials of $(n-2)^{\frac{1}{3}}$, so we can drop it out and take (\ref{48})
as the contribution of $u_{i}$, $1\leq i\leq r$. By (\ref{50}),
\ben
\beta_{M_i}(2n-4,j)+\beta_{M_i}(2n-4,n-2+j)=2\beta_{M_i}(n-2,j),
\een
Therefore the contribution of $v_1$ and all leaves neighboring to $v_1$ is

\bea\label{52}
&&2^{r}\sum_{j=0}^{2n-5}\beta_{M_{r+1}}(2n-4,j+k)\prod_{i=1}^{r}\beta_{M_i}(n-2,j)\frac{\pd}{\pd v_{l_i}^{j}}F_{0}^{\hat{D}_n}\Big|_{\text{all}\hspace{0.1cm} v_{s}^{t}=0}\nonumber\\
&=& 2^{r}\sum_{j=0}^{2n-5}\beta_{M_{r+1}}(2n-4,j+k)\prod_{i=1}^{r}\sum_{a_i=0}^{n-1}\zeta_{n-2}^{ja_{i}}\prod_{m\in M_{i}}B_{m}(\frac{a_i}{n-2})\frac{\pd}{\pd v_{l_i}^{j}}F_{0}^{\hat{D}_n}\Big|_{\text{all}\hspace{0.1cm} v_{s}^{t}=0}\nonumber\\
&=& \sum_{a_1=0}^{n-1}\cdots\sum_{a_r=0}^{n-1}\Bigg[\sum_{j=0}^{2n-5}\beta_{M_{r+1}}(2n-4,j+k)\prod_{i=1}^{r}\zeta_{n-2}^{ja_{i}}\frac{\pd}{\pd v_{l_i}^{j}}F_{0}^{\hat{D}_n}\Big|_{\text{all}\hspace{0.1cm} v_{s}^{t}=0}\cdot 2^{r}\prod_{p=1}^{r}\prod_{m\in M_{p}}B_{m}(\frac{a_p}{n-2})\Bigg].\nonumber\\
\eea
In the last equality we have separated the product in the squared brackets into two parts. We compute the first part individually, since this computation will be common in the induction.
\bea\label{55}
&&\sum_{j=0}^{2n-5}\beta_{M_{r+1}}(2n-4,j+k)\prod_{i=1}^{r}\Big(
\zeta_{n-2}^{ja_{i}}\frac{\pd}{\pd v_{l_i}^{j}}\Big)F_{0}^{\hat{D}_n}\Big|_{\text{all}\hspace{0.1cm} v_{s}^{t}=0}\nonumber\\
&=&\beta_{M_{r+1}}(2n-4,k)\prod_{i=1}^{r}(4n-8)^{\frac{2\sum l_i-2r}{3}}
\Bigg[(\pd_{l_1,\cdots,l_r}F_{0}^{\hat{D}_{n}})\Big(\Big\{\frac{(n-2)t_{l}^{\lbr b\rbr}+(n-2)t_{l}^{\lbr ab\rbr}
}{(4n-8)^{\frac{2(1-l)}{3}}}\Big\}\Big)\nonumber\\
&&+(\pd_{l_1,\cdots,l_r}F_{0}^{\hat{D}_{n}})\Big(\Big\{\frac{-(n-2)t_{l}^{\lbr b\rbr}-(n-2)t_{l}^{\lbr ab\rbr}
}{(4n-8)^{\frac{2(1-l)}{3}}}\Big\}\Big)\Bigg]\nonumber\\
&&+\beta_{M_{r+1}}(2n-4,k+n-2)\prod_{i=1}^{r}(4n-8)^{\frac{2\sum l_i-2r}{3}}\Bigg[(\pd_{l_1,\cdots,l_r}F_{0}^{\hat{D}_{n}})\Big(\Big\{\frac{(n-2)t_{l}^{\lbr b\rbr}+(n-2)t_{l}^{\lbr ab\rbr}
}{(4n-8)^{\frac{2(1-l)}{3}}}\Big\}\Big)\nonumber\\
&&+(\pd_{l_1,\cdots,l_r}F_{0}^{\hat{D}_{n}})\Big(\Big\{\frac{-(n-2)t_{l}^{\lbr b\rbr}-(n-2)t_{l}^{\lbr ab\rbr}
}{(4n-8)^{\frac{2(1-l)}{3}}}\Big\}\Big)\Bigg]\nonumber\\
&&+\sum_{j=1}^{n-3}\big(\beta_{M_{r+1}}(2n-4,j+k)+\beta_{M_{r+1}}(2n-4,n-2+j+k)\big)\nonumber\\
&&\cdot\prod_{i=1}^{r}
\zeta_{n-2}^{ja_{i}}(2n-4)^{\frac{2\sum l_i-2r}{3}}\pd_{l_1,\cdots,l_r}F_{0}^{\hat{D}_{n}}({0})\nonumber\\
&=&2\beta_{M_{r+1}}(n-2,k)\prod_{i=1}^{r}(4n-8)^{\frac{2\sum l_i-2r}{3}}\cdot\Bigg[(\pd_{l_1,\cdots,l_r}F_{0}^{\hat{D}_{n}})\Big(\Big\{\frac{(n-2)t_{l}^{\lbr b\rbr}+(n-2)t_{l}^{\lbr ab\rbr}
}{(4n-8)^{\frac{2(1-l)}{3}}}\Big\}\Big)\nonumber\\
&&+(\pd_{l_1,\cdots,l_r}F_{0})\Big(\Big\{\frac{-(n-2)t_{l}^{\lbr b\rbr}-(n-2)t_{l}^{\lbr ab\rbr}
}{(4n-8)^{\frac{2(1-l)}{3}}}\Big\}\Big)\Bigg]\nonumber\\
&&+\sum_{j=0}^{n-3}\beta_{M_{r+1}}(n-2,j+k)\prod_{i=1}^{r}
\zeta_{n-2}^{ja_{i}}(2n-4)^{\frac{2\sum l_i-2r}{3}}\pd_{l_1,\cdots,l_r}F_{0}^{\hat{D}_{n}}({0})\nonumber\\
&&-2\beta_{M_{r+1}}(n-2,k)\prod_{i=1}^{r}
(2n-4)^{\frac{2\sum l_i-2r}{3}}\pd_{l_1,\cdots,l_r}F_{0}^{\hat{D}_{n}}({0}).
\eea
Now we call back the second part of (\ref{52}) and compute the sum
\bea\label{53}
&&\sum_{a_1=0}^{n-1}\cdots\sum_{a_r=0}^{n-1}\sum_{j=0}^{n-3}\Big(\beta_{M_{r+1}}(n-2,j+k)\prod_{i=1}^{r}
\zeta_{n-2}^{ja_{i}}\cdot 2^{r}\prod_{m\in M_{i}}B_{m}(\frac{a_i}{n-2})\Big)\nonumber\\
&=& 2^{r}(n-2)\sum_{0\leq a_{i}\leq n-1, 1\leq i\leq r+1, a_{1}+\cdots a_{r+1}\equiv 0(\text{mod}\hspace{0.1cm}n)} \zeta_{n-2}^{ka_{r+1}}
\prod_{i=1}^{r+1}\prod_{m\in M_{i}}B_{m}(\frac{a_i}{n-2}),\nonumber\\
\eea

Also,
\bea\label{54}
&&\sum_{a_1=0}^{n-1}\cdots\sum_{a_r=0}^{n-1}\Big(\beta_{M_{r+1}}(n-2,k)\prod_{i=1}^{r}
 2^{r}\prod_{m\in M_{i}}B_{m}(\frac{a_i}{n-2})\Big)\nonumber\\
&=& 2^{r}\beta_{M_{r+1}}(n-2,k)\beta_{M_{1}}(n-2,0)\cdots \beta_{M_{r}}(n-2,0)\nonumber\\
&=& 2^{r}\beta_{M_{1}}(n-2,0)\cdots \beta_{M_{r}}(n-2,0)\cdot\sum_{a_{r+1}}\zeta_{n}^{ka_{r+1}}\prod_{m\in M_{r+1}}B_{m}(\frac{a_{r+1}}{n-2}),\nonumber\\
\eea
by lemma \ref{22}, $\beta_{M_{1}}(n-2,0)\cdots \beta_{M_{r}}(n-2,0)$ is a Laurent polynomial of $n-2$. By (\ref{55}), (\ref{53}),
(\ref{54}), the \emph{output} of the contribution of $v_1$ is a linear combination of
\bea
\sum_{0\leq a_{i}\leq n-1, 1\leq i\leq r+1, a_{1}+\cdots a_{r+1}\equiv 0(\text{mod}\hspace{0.1cm}n)} \zeta_{n-2}^{ka_{r+1}}
\prod_{i=1}^{r+1}\prod_{m\in M_{i}}B_{m}(\frac{a_i}{n-2})
\eea
and
\bea
\sum_{a_{r+1}}\zeta_{n}^{ka_{r+1}}\prod_{m\in M_{r+1}}B_{m}(\frac{a_{r+1}}{n-2}),
\eea
with coefficients Laurent polynomials of $(n-2)^{\frac{1}{3}}$. Then one does the same thing for all vertices in $V(\Gamma(1))\backslash V(\Gamma(2))$. In the induction from $\Gamma^{(1)}$ to $\Gamma^{(2)}$, the same computations as (\ref{55}) applies. And one finally sees that we obtain many summing as $S_{n-2}(\Gamma_i)$ for subtrees $\Gamma_{i}$. \\

As a last word, one need not worry about the fractional power of $(n-1)^{\frac{1}{3}}$, since the fractional power in (\ref{56}) arises as an integer power distributed to several variables, by the proof of proposition 4.2 in \cite{JK}. Thus we have proved the Claim 2. \hfill\qedsymbol\\

Finally we need consider the operators of the form (\ref{40}) and (\ref{41}) together. However, including the operators of the form (\ref{40}) just gives more times of differentiation to  $F_{0}^{\hat{D}_n}$ in (\ref{55}), and one easily sees that the Claim 2 still holds. Therefore we complete the proof of the theorem.\hfill\qedsymbol

\begin{remark}
In fact our proof gives more than we need, i.e., we have proved the polynomiality of the \emph{gravitational descendant} correlators $\langle \prod_{i=1}^{k}\tau_{p_i}(e_{\lbr b \rbr}) \prod_{j=1}^{l}\tau_{q_j}(e_{\lbr ab \rbr})\rangle^{[\mathbb{C}^2 /\hat{D}_n]}$ for fixed $k, l$ and $p_i, q_j$ independent of $n$. This result coincides with the theorem 2.7 in \cite{Maulik} for suitable \emph{reduced} descendant invariants, by the corresponding conjectural quantum McKay correspondence for descendant invariants.
\end{remark}

\begin{remark}
One can use the same method to prove the polynomiality of $\langle e_{\lbr b \rbr}^{2m} \rangle^{[\mathbb{C}^2 /\hat{D}_n]}$ for odd $n$. In fact  the proof for odd $n$ is less complicated, as in the WDVV induction. We only treat the cases for even $n$ because in corollary \ref{30} we showed that the quantum McKay correspondence holds for infinitely many even $n$.
\end{remark}

\hspace{1cm}\footnotesize{Department of Mathematical Sciences, Tsinghua University, Beijing, 100084, China }\\

\hspace{1cm}\footnotesize{\emph{E-mail address}: huxw08@mails.tsinghua.edu.cn}
\end{document}